\documentclass[10pt]{article}

\usepackage{amsmath,amssymb,amsfonts,amsthm,epic}
\usepackage[vcentermath]{youngtab}
\usepackage{pstricks,pst-node}

\theoremstyle{plain}
\newtheorem{thm}[subsection]{Theorem}
\newtheorem{lem}[subsection]{Lemma}
\newtheorem{prop}[subsection]{Proposition}
\newtheorem{cor}[subsection]{Corollary}
\newtheorem{ex}[subsection]{Example}

\theoremstyle{definition}

\newcommand{\es}{\emptyset}
\newcommand{\bN}{\mathbb{N}}
\newcommand{\bC}{\mathbb{C}}

\newcommand{\ga}{\alpha}
\newcommand{\gb}{\beta}
\newcommand{\gc}{\gamma}
\newcommand{\gi}{\iota}

\newcommand{\gs}{\sigma}
\newcommand{\gt}{\theta}

\newcommand{\cO}{\mathcal{O}}

\newcommand{\cT}{\mathcal{T}}

\newcommand{\cM}{\mathcal{M}}
\newcommand{\disjunion}{\,\,\dot{\cup}\,\,}

\newcommand{\NN}{\bN^2}

\newcommand{\gB}{\ol{\gb}\times\gb}
\newcommand{\BBZ}{\gB}
\newcommand{\BBP}{(\gB)^+}
\newcommand{\BBN}{(\gB)^-}

\newcommand{\pr}{R}
\newcommand{\qr}{S}

\newcommand{\rsk}{\mathop{\rm BRSK}\nolimits}
\newcommand{\rrsk}{\mathop{\rm RBRSK}\nolimits}

\newcommand{\geqT}{\geq}
\newcommand{\leqT}{\leq}
\newcommand{\gtT}{>}
\newcommand{\ltT}{<}
\newcommand{\gtS}{\gtrdot}
\newcommand{\ltS}{\lessdot}

\newcommand{\ltC}{\prec}

\newcommand{\geqCC}{\unrhd}
\newcommand{\leqCC}{\unlhd}

\newcommand{\ol}{\overline}

\newcommand{\wt}{\widetilde}

\newcommand{\la}[1]{\stackrel{#1}{\leftarrow}}
\newcommand{\ra}[1]{\stackrel{#1}{\rightarrow}}

\newcommand{\floor}[1]{\lfloor #1 \rfloor}
\newcommand{\ceil}[1]{\lceil #1 \rceil}
\newcommand{\sfloor}[1]{\lfloor\! #1 \!\rfloor}
\newcommand{\sceil}[1]{\lceil\! #1 \!\rceil}

\newcommand{\meet}{\wedge}

\newcommand{\Span}{\mathop{\rm Span}\nolimits}
\newcommand{\depth}{\mathop{\rm depth}\nolimits}

\newcommand{\Mult}{\mathop{\rm Mult}\nolimits}

\newcommand{\abc}{_{\ga,\gb}^\gc}

\newcommand{\init}{\mathop{\rm in}\nolimits}
\newcommand{\pl}{\mathop{\rm pl}\nolimits}

\newcommand{\plotstarchar}{\makebox(0,0){\large $\times$}}
\newcommand{\plotcirchar}{\makebox(2,-1){\circle*{2}}}
\newcommand{\plotsmcirchar}{\makebox(.8,-.717){\circle*{.5}}}
\newcommand{\plotsmbcirchar}{\makebox(1,-1){$\blacksquare $}}

\newcommand{\bsq}{{\bullet}}
\newcommand{\pone}{{\phantom{1}}}

\title{Local Properties of Richardson Varieties in the
Grassmannian via a Bounded Robinson-Schensted-Knuth
Correspondence}
\date{}
\author {Victor Kreiman}

\begin{document}
\maketitle

\begin{abstract}
The Richardson variety $X_\ga^\gc$ in the Grassmannian is defined
to be the intersection of the Schubert variety $X^\gc$ and
opposite Schubert variety $X_\ga$.  We give an explicit Gr\"obner
basis for the ideal of the tangent cone at any $T$-fixed point of
$X_\ga^\gc$, thus generalizing a result of Kodiyalam-Raghavan
\cite{Ko-Ra} and Kreiman-Lakshmibai \cite{Kr-La}. Our proof is
based on a generalization of the Robinson-Schensted-Knuth (RSK)
correspondence, which we call the bounded RSK (BRSK).  We use the
Gr\"obner basis result to deduce a formula which computes the
multiplicity of $X_\ga^\gc$ at any $T$-fixed point by counting
families of nonintersecting lattice paths, thus generalizing a
result first proved by Krattehthaler \cite{Kra1,Kra2}.
\end{abstract}

\setcounter{tocdepth}{1} \tableofcontents

\section{Introduction}
The Richardson variety $X_\ga^\gc$ in the
Grassmannian\footnote{Richardson varieties in the Grassmannian are
also studied by Stanley in \cite{St}, where these varieties are
called {\it skew Schubert varieties}.  Discussion of these
varieties also appears in \cite{Ho-Pe}.} is defined to be the
intersection of the Schubert variety $X^\gc$ and opposite Schubert
variety $X_\ga$.  In particular, Schubert and opposite Schubert
varieties are special cases of Richardson varieties.  We derive
local properties of $X_\ga^\gc$ at any $T$-fixed point $e_\gb$. It
should be noted that the local properties of Schubert varieties at
$T$-fixed points determine their local properties at all other
points, because of the $B$-action; but this does not extend to
Richardson varieties, since Richardson varieties only have a
$T$-action.

In Kodiyalam-Raghavan \cite{Ko-Ra} and Kreiman-Lakshmibai
\cite{Kr-La}, an explicit Gr\"obner basis for the ideal of the
tangent cone of the Schubert variety $X^\gc$ at $e_\gb$ is
obtained. The Gr\"obner basis is used to derive a formula for the
multiplicity of $X^\gc$ at $e_\gb$. In this paper, we generalize
the results of \cite{Ko-Ra} and \cite{Kr-La} to the case of
Richardson varieties.
The results of \cite{Ko-Ra} and \cite{Kr-La} were conjectured by
Kreiman-Lakshmibai \cite{Kr-La2}, although in a different,
group-theoretic form. The multiplicity formula (in both forms) was
first proved by Krattenthaler \cite{Kra1, Kra2} by
showing its equivalence to the Rosenthal-Zelevinsky determinantal
multiplicity formula \cite{Ro-Ze}.


Sturmfels \cite{Stu} and Herzog-Trung \cite{He-Tr} proved results
on a class of determinantal varieties which are equivalent to the
results of \cite{Ko-Ra}, \cite{Kr-La}, and this paper for the case
of Schubert varieties at the $T$-fixed point $e_{id}$. The key to
their proofs was to use a version of the Robinson-Schensted-Knuth
correspondence
(which we shall call the `ordinary' RSK)
in order to establish a degree-preserving bijection
between a set of monomials defined by an initial ideal and a
`standard monomial basis'.

The difficulty in generalizing this method of proof to the case of
Schubert varieties at an arbitrary $T$-fixed point $e_\gb$ lies in
generalizing this bijection. All three of \cite{Ko-Ra},
\cite{Kr-La}, and this paper obtain generalizations of this
bijection.
The three generalizations, when restricted to Schubert varieties,
are in fact the same bijection\footnote{This supports the
conviction of the authors in \cite{Ko-Ra} that this bijection is
natural and that it is in some sense the only natural bijection
satisfying the required geometric conditions.}, although this is
not immediately apparent. Although the formulations of the
bijection in \cite{Ko-Ra} and \cite{Kr-La} are similar to
eachother, our formulation is in terms of different combinatorial
indexing sets, and thus most of our combinatorial definitions and
proofs are of a different nature than those of \cite{Ko-Ra} and
\cite{Kr-La}. The relationship between our formulation and the
formulations in \cite{Ko-Ra} and \cite{Kr-La} is analogous to the
relationship between the Robinson-Schensted correspondence and
Viennot's version of the Robinson-Schensted correspondence
\cite{Sa,Vi}.

We formulate the bijection by introducing a generalization of the
ordinary RSK, which we call the bounded RSK. Because the
definition of the bounded RSK is built from that of the ordinary
RSK, many properties of the bounded RSK are immediate consequences
of analogous properties for the ordinary RSK (see \cite{Fu, Sa}).
This simplifies our proofs.

Results analogous to those of \cite{Ko-Ra} and \cite{Kr-La} have
now been obtained for the symplectic and orthogonal Grassmannians
(see \cite{Gh-Ra}, \cite{Ra-Up}). We believe it is possible that
the methods of this paper can be adapted to these varieties as
well. The results of \cite{Ko-Ra}, \cite{Kr-La}, and this paper
have been used to study the equivariant cohomology and equivariant
K-theory of the Grassmannian (see \cite{Kr}, \cite{La-Ra-Sa}).

\vspace{1em}

\noindent\textbf{Acknowledgements.} I would like to thank A.
Conca, V. Lakshmibai, P. Magyar, B. Sagan, and M. Shimozono for
valuable discussions, suggestions, and corrections. I am also
grateful to the referees for their careful readings and for
pointing out a number of improvements and corrections.

\section{Statement of Results}\label{s.results}

Let $K$ be an algebraically closed field, and let $d$, $n$ be
fixed positive integers, $0<d<n$. The \textbf{Grassmannian}
$Gr_{d,n}$ is the set of all $d$-dimensional subspaces of $K^n$.
The \textbf{Pl\"ucker map} $\pl:Gr_{d,n}\to\mathbb{P}(\wedge^d
K^n)$ is given by $\pl(W)=[w_1\wedge\cdots\wedge w_d]$, where
$\{w_1,\ldots,w_d\}$ is any basis for $W$. It is well known that
$\pl$ is a bijection onto a closed subset of $\mathbb{P}(\wedge^d
K^n)$. Thus $Gr_{d,n}$ inherits the structure of a projective
variety.

Define $I_{d,n}$ to be the set of $d$-element subsets of
$\{1,\ldots,n\}$. Let $\ga=\{\ga_1,\ldots,\ga_d\}\in I_{d,n}$,
$\ga_1<\cdots<\ga_d$. Define the \textbf{complement} of $\ga$ by
$\ol{\ga}=\{1,\ldots,n\}\setminus\ga$, and the \textbf{length} of
$\ga$ by $l(\ga)=\ga_1+\cdots+\ga_d - {d+1\choose 2}$. If
$\gb=\{\gb_1,\ldots,\gb_d\}\in I_{d,n}$, $\gb_1<\ldots<\gb_d$,
then we say that $\ga\leq \gb$ if $\ga_i\leq\gb_i$,
$i=1,\ldots,d$.

Let $\{e_1,\ldots,e_n\}$ be the standard basis for $K^n$. For
$\ga\in I_{d,n}$, define
$e_\ga=\Span\{e_{\ga_1},\ldots,e_{\ga_d}\}\in Gr_{d,n}$.
Let
\begin{align*}
&G=GL_n(K)\\
&B=\{g\in G\mid g\hbox{ is upper triangular}\}\\
&B^-=\{g\in G\mid g\hbox{ is lower triangular}\}\\
&T=\{g\in G\mid g\hbox{ is diagonal}\}
\end{align*}
The group $G$ acts transitively on $Gr_{d,n}$ with $T$-fixed
points $\{e_\ga\mid \ga\in I_{d,n}\}$. The Zariski closure of the
$B$ (resp. $B^-$) orbit through $e_\ga$, with canonical reduced
scheme structure, is called a \textbf{Schubert variety} (resp.
\textbf{opposite Schubert variety}), and denoted by $X^\ga$ (resp.
$X_\ga$).  For $\ga,\gc\in I_{d,n}$, the scheme-theoretic
intersection $X_\ga^\gc=X_\ga\cap X^\gc$ is called a
\textbf{Richardson variety}. It can be shown that $X_\ga^\gc$ is
nonempty if and only if $\ga\leq \gc$; that $e_\gb\in X_\ga^\gc$
if and only if $\ga\leq\gb\leq\gc$; and that if $X_\ga^\gc$ is
nonempty, it is reduced and irreducible of dimension
$l(\gc)-l(\ga)$ (see \cite{Br,Kr-La3,La-Go,Ri}).

For $\gb\in I_{d,n}$ define $p_\gb$ to be homogeneous
(\textbf{Pl\"ucker}) coordinate $[e_{\gb_1}\wedge\cdots\wedge
e_{\gb_d}]^*$ of $\mathbb{P}(\wedge^d K^n)$. Let $\cO_\gb$ be the
distinguished open set of $Gr_{d,n}$ defined by $p_\gb\neq 0$. Its
coordinate ring $K[\cO_\gb]$ is isomorphic to the homogeneous
localization $K[Gr_{d,n}]_{(p_{\gb})}$. Define $f_{\gt,\gb}$ to be
$p_\gt/p_\gb\in K[\cO_\gb]$.

The open set $\cO_\gb$ is isomorphic to the affine space
$K^{d(n-d)}$. Indeed, it can be identified with the space of
matrices in $M_{n\times d}$ in which rows $\gb_1,\ldots,\gb_d$ are
the rows of the $d\times d$ identity matrix, and rows
$\ol{\gb}_1,\ldots,\ol{\gb}_{n-d}$ contain arbitrary elements of
$K$. The rows of $\cO_\gb$ are indexed by $\{1,\ldots,n\}$, and
the columns by $\gb$.  Note that the affine coordinates of
$K[\cO_\gb]$ are indexed by $\BBZ$.  The coordinate $f_{\gt,\gb}$,
$\gt\in I_{d,n}$, is identified with plus or minus the $d\times d$
minor of $\cO_\gb$ with row-set $\gt_1,\ldots,\gt_d$.

\begin{ex}\label{ex.big-cell}
Let $d=3$, $n=7$, $\gb=\{2,5,7\}$. Then
\begin{align*}
\cO_\gb&=\left\{ \left(
\begin{matrix}
x_{12}&x_{15}&x_{17}\\
1&0&0\\
x_{32}&x_{35}&x_{37}\\
x_{42}&x_{45}&x_{47}\\
0&1&0\\
x_{62}&x_{65}&x_{67}\\
0&0&1
\end{matrix}
\right), x_{ij}\in K \right\},\\
\intertext{and}
f_{\{1,4,5\},\gb}&=\left|\begin{matrix}x_{12}&x_{15}&x_{17}\\
x_{42}&x_{45}&x_{47}\\ 0&1&0\end{matrix}\right|=-\left|\begin{matrix}x_{12}&x_{17}\\
x_{42}&x_{47}\end{matrix}\right|.
\end{align*}
\end{ex}

\vspace{2em}

In order to better understand the local properties of $X_\ga^\gc$
near $e_\gb$, we analyze $Y_{\ga,\gb}^\gc=X_\ga^\gc\cap \cO_\gb$,
an open subset of $X_{\ga}^\gc$ centered at $e_\gb$, and a closed
affine subvariety of $\cO_\gb$. Let
$G_{\ga,\gb}^\gc=\{f_{\gt,\gb}\mid \ga\not\geqT \gt\hbox{ or }
\gt\not\leqT \gc\}\subset K[\cO_\gb]$, and let $\langle
G\abc\rangle$ be the ideal of $K[\cO_\gb]$ generated by $G\abc$.
The following is a well known result (see \cite{Br,La-Go}, for
instance).
\begin{thm}\label{t.r.ideal_gens} $K[Y_{\ga,\gb}^\gc]=K[\cO_\gb]/\langle
G\abc\rangle$.
\end{thm}

As a consequence of Theorem \ref{t.r.ideal_gens},
$Y_{\ga,\gb}^\gc$ is isomorphic to the tangent cone of
$X_{\ga}^\gc$ at $e_\gb$, and thus $\deg
Y_{\ga,\gb}^\gc=\Mult_{e_\gb} X_\ga^\gc$, the multiplicity of
$X_\ga^\gc$ at $e_\gb$. Indeed, since $Y_{\ga,\gb}^\gc$ is an
affine variety in $\cO_\gb$ defined by a homogeneous ideal, with
$e_\gb$ the origin of $\cO_\gb$, $Y_{\ga,\gb}^\gc$ is isomorphic
to the tangent cone of $Y_{\ga,\gb}^\gc$ at $e_\gb$; since
$Y_{\ga,\gb}^\gc$ is open in $X_{\ga}^{\gc}$, the tangent cone of
$Y_{\ga,\gb}^\gc$ at $e_\gb$ is isomorphic to the tangent cone of
$X_{\ga}^{\gc}$ at $e_\gb$.

Any minor $f_{\gt,\gb}$ can be expressed naturally as plus or
minus a determinant all of whose entries are $x_{ij}$'s.  Choose a
monomial order on $K[\cO_\gb]$ such that the initial term of any
minor $f_{\gt,\gb}$, $\init(f_{\gt,\gb})$, is the
Southwest-Northeast diagonal of this determinant. The main result
of this paper, which is also proven in \cite{Ko-Ra} and
\cite{Kr-La}, is the following.
\begin{prop}\label{p.i.grobner}
$G_{\ga,\gb}^\gc$ is a Gr\"obner basis for $\langle G\abc\rangle$.
\end{prop}

If $S$ is any subset of $K[\cO_\gb]$, define $\init S$ to be the
ideal $\langle\init(s)\mid s\in S\rangle$.
\begin{cor}\label{c.i.grobner_extras} $\deg Y\abc$ ($=\Mult_{e_\gb}X_\ga^\gc$)  is the number of
square-free monomials of degree $l(\gc)-l(\ga)$ in
$K[\cO_\gb]\setminus \init G\abc$.
\end{cor}

We now briefly sketch the proof of Proposition \ref{p.i.grobner}
(omitting the details, which can be found in Section
\ref{s.grobner}), in order to introduce the main combinatorial
objects of interest and outline the structure of this paper. We
wish to show that in any degree, the number of monomials of $\init
G\abc$ is at least as great as the number of monomials of $\init
\langle G\abc\rangle$ (the other inequality being trivial), or
equivalently, that in any degree, the number of monomials of
$K[\cO_\gb]\setminus\init G\abc$ is no greater than the number of
monomials of $K[\cO_\gb]\setminus\init \langle G\abc\rangle$. Both
the monomials of $K[\cO_\gb]\setminus \init \langle G\abc\rangle$
and the standard monomials on $Y\abc$ form a basis for
$K[\cO_\gb]/ \langle G\abc\rangle$, and thus agree in cardinality
in any degree. Therefore, it suffices to give a degree-preserving
injection from the monomials of $K[\cO_\gb]\setminus\init G\abc$
to the standard monomials on $Y\abc$.  We construct such an
injection, the bounded RSK ($\rsk$), from an indexing set of the
former to an indexing set of the
latter. These indexing sets are given in the following table.\\

\begin{figure}[!h]
\begin{center}
\begin{tabular}{|l|l|}\hline
\multicolumn{1}{|c|}{\slshape Set of elements in $K[\cO_\gb]$}&
\multicolumn{1}{c|}{\slshape Indexing set}\\
\hline\hline
&{nonvanishing multisets}\\
\raisebox{.5em}[0pt]{monomials of $K[\cO_\gb]\setminus \init G\abc$}& on $\BBZ$ bounded by $T_\ga, W_\gc$\\
\hline
&{nonvanishing semistandard notched bitableaux}\\
\raisebox{.5em}[0pt]{standard monomials on $Y\abc$}&on $\BBZ$
bounded by $T_\ga, W_\gc$
\\
\hline
\end{tabular}
\end{center}
\caption{\label{f.i.main_objects}Two subsets of $K[\cO_\gb]$ and
their indexing sets}
\end{figure}

\noindent In Sections \ref{s.notched_tableaux}, \ref{s.multisets},
\ref{s.ss_notched_bitab}, \ref{s.grsk}, and \ref{s.grsk_props}, we
define nonvanishing multisets on $\BBZ$ bounded by $T_\ga, W_\gc$,
nonvanishing semistandard notched bitableaux on $\BBZ$ bounded by
$T_\ga, W_\gc$, and the injection $\rsk$ from the former to the
latter. In Section \ref{s.grobner}, we prove that these two
combinatorial objects are indeed indexing sets for the monomials
of $K[\cO_\gb]\setminus \init G\abc$ and the standard monomials on
$Y\abc$ respectively, and use this to prove Proposition
\ref{p.i.grobner} and Corollary \ref{c.i.grobner_extras}. In
Sections \ref{s.twisted} and \ref{s.path_families_mults}, we show
how using Corollary \ref{c.i.grobner_extras},
$\Mult_{e_\gb}X_\ga^\gc$ can be interpreted as counting certain
families of nonintersecting paths in the lattice $\gB$. In Section
\ref{s.proofs}, we give two of the more detailed proofs.

\section{Notched Tableaux and Bounded Insertion}\label{s.notched_tableaux}

A \textbf{Young diagram} (resp. \textbf{notched diagram}) is a
collection of boxes arranged into a left and top justified array
(resp. into left justified rows).  The \textbf{empty Young
diagram} is the Young diagram with no boxes.  A notched diagram
may contain rows with no boxes; however, a Young diagram may not,
unless it is the empty Young diagram. A \textbf{Young tableau}
(resp. \textbf{notched tableau}) is a filling of the boxes of a
Young diagram (resp. notched diagram) with positive integers. The
\textbf{empty Young tableau} is the Young tableau with no boxes.
Let $P$ be either a notched tableau or a Young tableau. We denote
by $P_i$ the $i$-th row of $P$ from the top, and by $P_{i,j}$ the
$j$-th entry from the left of $P_i$. We say that $P$ is
\textbf{row strict} if the entries of any row of $P$ strictly
increase as you move to the right.  If $P$ is a Young tableau,
then we say that $P$ is \textbf{semistandard} if it is row strict
and the entries of any column weakly increase as you move down. By
definition, the empty Young tableau is considered semistandard.
\begin{ex}
A row strict notched tableau $P$, and a semistandard Young tableau
$R$.
\[
P=\young(2346,78,1679,689,34567,234579)\ , \qquad
R=\young(12456,13478,235,26,27,9)\\[1em]
\]
\end{ex}

Let $P$ be a row strict notched tableau, and $b$ a positive
integer. Since $P$ is row strict, its entries which are greater
than or equal to $b$ are right justified in each row.  Thus if we
remove these entries (and their boxes) from $P$ then we are left
with a row strict notched tableau, which we denote by $P^{<b}$. We
say that $P$ is \textbf{semistandard on} $\mathbf{b}$ if $P^{<b}$
is a semistandard Young tableau (note that if $P$ is semistandard
on $b$ and the first row of $P^{<b}$ has no boxes, then $P^{<b}$
must be the empty Young tableau). It is clear that if $P$ is
semistandard on $b$, then it is semistandard on $b'$ for any
positive integer $b'<b$.
\begin{ex}
Let $P=\young(1246,236,24578,3,45)$.  Then
$$P^{<5}=\young(124,23,24,3,4),\qquad
P^{<6}=\young(124,23,245,3,45)\ .$$ Thus $P$ is semistandard on
$5$ but not on 6.

\end{ex}

We next review the transpose of Schensted's column insertion
process, which we shall call simply `ordinary' Schensted
insertion. It is an algorithm which takes as input a semistandard
Young tableau $P$ and a positive integer $a$, and produces as
output a new semistandard Young tableau with the same shape as $P$
plus one extra box, and with the same entries as $P$ (possibly in
different locations) plus one additional entry, namely $a$. To
begin, insert $a$ into the first row of $R$, as follows. If $a$ is
greater than all entries in the first row of $R$, then place $a$
in a new box on the right end of the first row, and ordinary
Schensted insertion terminates. Otherwise, find the smallest entry
of the first row of $R$ which is greater than or equal to $a$, and
replace that number with $a$. We say that the number which was
replaced was ``bumped'' from the first row. Insert the bumped
number into the second row in precisely the same fashion that $a$
was inserted into the first row. This process continues down the
rows until, at some point, a number is placed in a new box on the
right end of some row, at which point ordinary Schensted insertion
terminates.

We next describe the \textbf{bounded insertion algorithm}, which
takes as input a positive integer $b$, a notched tableau $P$ which
is semistandard on $b$, and a positive integer $a<b$, and produces
as output a notched tableau which is semistandard on $b$, which we
denote by $P\la{b}a$.

\begin{quote}
\noindent\textbf{Bounded Insertion} {\it
\begin{itemize}

\item[\textbf{Step 1.}] Remove all entries of $P$ which are
greater than or equal to $b$ from $P$, resulting in the
semistandard Young tableau $P^{<b}$.

\item[\textbf{Step 2.}] Insert $a$ into $P^{<b}$ using the
ordinary Schensted insertion process.

\item[\textbf{Step 3.}] Place the entries of $P$ which were
removed when forming $P^{<b}$ in Step 1 back into the Young
tableau resulting from Step 2, in the same rows from which they
were removed.
\end{itemize}
}
\end{quote}
This insertion process is effectively the ordinary Schensted
insertion of $a$ into $P$, but acting only on the part of $P$
which is ``bounded'' by $b$. The fact that bounded insertion
preserves the property of being semistandard on $b$ follows
immediately from the fact that ordinary Schensted insertion
preserves the
property of being semistandard.\\

\begin{ex}\label{e.insertion}
Let $P=\young(1247,158,36789,46)$, $a=3$, $b=6$. We compute
$P\la{b}a$.\\[1em]

\noindent Step 1. Remove all entries of $P$ which are greater
than or equal than $b$, resulting in $$P^{<b}=\young(124,15,3,4).$$\\

\noindent Step 2. Insert $a$ into $P^{<b}$ using the ordinary
Schensted insertion process: $a$ bumps 4 from the first row, which
bumps 5 from the second row, which is placed in a new box on the
right end of the third row, to form
$$\young(123,14,35,4).$$\\

\noindent Step 3. Place the entries removed from $P$ in Step 1
back into the Young tableau resulting from Step 2, in the same
rows from which they were removed, to obtain
$$P\la{b}a=\young(1237,148,356789,46)\ .$$
\end{ex}
\vspace{1em}

We define the \textbf{bumping route} of the bounded insertion
algorithm to be the sequence of boxes in $P$ from which entries
are bumped in Step 2, together with the \textbf{new box} which is
added at the end of Step 2.

\begin{ex}
The bumping route in Example \ref{e.insertion} is the set of boxes
with $\bsq$'s in the following Young diagram:
$$\young(\pone\pone\bsq\pone,\pone\bsq\pone,\pone\bsq\pone\pone\pone\pone,\pone\pone)$$
\vspace{1em} The new box is the lowest box containing a $\bsq$.
\end{ex}

The bounded insertion algorithm is reversible: if $P\la{b}a$ is
computed, and we know $b$ and the location of the new box, then we
can retrieve $P$ and $a$ by running the bounded insertion
algorithm in reverse. Note that the reverse of Step 2 is the
ordinary Schensted reverse insertion process.

Suppose that $P$ is semistandard on $b$, that $a,a'<b$, and that
bounded insertion is performed twice in succession, resulting in
$(P\la{b}a) \la{b}a'$. Let $R$ and $B$ be the bumping route and
new box of the first insertion, and let $R'$ and $B'$ be the
bumping route and new box of the second insertion. We say that
$R'$ is \textbf{weakly left} of $R$ if for every box of $R$, there
is a box of $R'$ to the left of or equal to it; we say that $R$ is
\textbf{strictly left} of $R'$ if for every box of $R'$, there is
a box of $R$ to the left of it. We say that $B'$ is
\textbf{strictly below} $B$ if $B'$ lies in a lower row than $B$;
we say that $B$ is \textbf{weakly below} $B'$ if $B$ lies in
either the same row as $B'$ or a lower row than $B'$. The
following Lemma is an immediate consequence of the analogous
result for ordinary Schensted insertion (see \cite{Fu}).

\begin{lem}\label{l.bumping_route}
(i) If $a\geq a'$, then $R'$ is weakly left of $R$ and $B'$ is
strictly below $B$.\\
(ii) If $a<a'$, then $R$ is strictly left of $R'$ and $B$ is
weakly below $B'$.
\end{lem}

\section{Multisets on $\bN$ and on $\bN^2$}\label{s.multisets}

Let $S$ be any set.  A multiset $E$ on $S$ is defined to be a
function $E:S\to \{0,1,2,\ldots,\}$.  One should think of $E$ as
consisting of the set $S$ of elements, but with each $s\in S$
occurring $E(s)$ times. Note that a set is a special type of
multiset in which each element occurs exactly once.  Define the
\textbf{underlying set of} $E$ to be $\{s\in S\mid E(s)\neq 0\}$,
a subset of $S$. If $T$ is a subset of $S$, then we write
$E\subset T$ if the underlying set of $E$ is a subset of $T$. We
often write a multiset $E$ by listing its elements,
$E=\{e_1,e_2,e_3,\ldots\}$, where the $e_i$'s may not be distinct
(in fact, each $e_i$ occurs $E(e_i)$ times in the list).

We call $E(s)$ the \textbf{degree} or \textbf{multiplicity} of $s$
in $E$. The multiset $E$ is said to be \textbf{finite} if $E(s)$
is nonzero for only finitely many $s\in S$. If $E$ is finite, then
define the \textbf{degree} of $E$, denoted by $|E|$, to be the sum
of $E(s)$ over all $s\in S$. If $E$ is not finite, then define the
degree of $E$ to be $\infty$.
Define the multisets $E\disjunion F$ and $E\setminus F$ as
follows:
\begin{align*}
(E\disjunion F)(s)&=E(s)+F(s),\ s\in S\\
(E\setminus F)(s)&=\max\{E(s)-F(s),0\},\ s\in S
\end{align*}
\begin{ex}
Let $E=\{a,b,b,b,c,c,c\}$, $F=\{b,b,c,d\}$.  Then $|E|=7$,
$E\disjunion F=\{a,b,b,b,b,b,c,c,c,c,d\}$, and $E\setminus
F=\{a,b,c,c\}$.
\end{ex}

\subsection*{Multisets on $\bN$}

Let $\bN$ denote the positive integers. Let $E=\{e_1,e_2,\ldots\}$
be a multiset on $\bN$, and let $z\in \bN$. Define the multiset
$E^{\leq z}:=\{e_j\in E\mid e_j\leq z\}$.

Let $A=\{a_1,a_2,\ldots\}$ and $B=\{b_1,b_2,\ldots\}$ be two
multisets on $\bN$ of the same degree, with $a_i\leq a_{i+1}$,
$b_i\leq b_{i+1}$, for all $i$. We say that $A$ is less than or
equal to $B$ in the \textbf{termwise order} if $a_i\leq b_i$ for
all $i$, or equivalently if $|A^{< z}|\geq |B^{< z}|$ for all
$z\in\bN$.  We denote this by $A\leq B$. We say that $A$ is less
than $B$ in the \textbf{strict termwise order} if $a_i< b_i$ for
all $i$.  We denote this by $A\ltS B$. Note that $\leq$ is a finer
order than $\ltS$.

If $A$, $B$, $C$, and $D$ are multisets on $\bN$ such that
$|A\disjunion D|=|B\disjunion C|$, then we write
\begin{equation}\label{e.m.set_subtraction}
A-C\leq B-D \hbox{ to indicate that }A\disjunion D\leq B\disjunion
C.
\end{equation}
Note that $A-B\leq C-D$ is a transitive relation.

In general no meaning is attached to the expression $A-C$ by
itself. However, if $A$ and $C$ are both sets, then we may define
$A-C:=A\disjunion (\bN\setminus C)$. If in addition
$A\subset\ol{\gb}$ and $C\subset\gb$, then we may define
$A-C:=A\disjunion(\gb\setminus C)$. It is an easy check that both
of these definitions are consistent with
(\ref{e.m.set_subtraction}) (e.g., $A\disjunion (\bN\setminus
C)\leq B\disjunion (\bN\setminus D)$ if and only if $A\disjunion
D\leq B\disjunion C$).

\subsection*{Multisets on $\bN^2$}

Let  $U=\{(e_1,f_1),(e_2,f_2),\ldots\}$ be a multiset on $\bN^2$.
Define $U_{(1)}$ and $U_{(2)}$ to be the multisets
$\{e_1,e_2,\ldots\}$ and $\{f_1,f_2,\ldots\}$ respectively on
$\bN$. Define the \textbf{nonvanishing}, \textbf{negative}, and
\textbf{positive parts} of $U$ to be the following multisets:
\begin{align*}
U^{\neq 0}&=\{(e_i,f_i)\in U\mid e_i-f_i\neq 0\},\\
U^-&=\{(e_i,f_i)\in U\mid e_i-f_i<0\},\\
U^+&=\{(e_i,f_i)\in U\mid e_i-f_i>0\}.
\end{align*}

\noindent It is useful to visualize the $e$-axis pointing downward
and the $f$-axis pointing to the right, as illustrated below (the
large squares cover the points of $\NN\setminus(\NN)^{\neq 0}$):

\begin{center}
\setlength{\unitlength}{.54mm}
\begin{picture}(90,110)(0,0)
\matrixput(10,20)(10,0){7}(0,10){7}{\plotsmcirchar}
\multiputlist(0,10)(0,10){$\vdots$,7,6,5,4,3,2,1}
\multiputlist(10,90)(10,0){1,2,3,4,5,6,7,$\cdots$}
\multiput(10,80)(10,-10){7}{\plotsmbcirchar}

\put(55,65){\makebox(0,0){$(\NN)^-$}}
\put(25,35){\makebox(0,0){$(\NN)^+$}}
\put(-15,50){\makebox(0,0){$e_i$'s}}
\put(40,105){\makebox(0,0){$f_i$'s}}

\end{picture}
\end{center}

We say that $U$ is \textbf{nonvanishing} if $U\subset (\NN)^{\neq
0}$, \textbf{negative} if $U\subset (\NN)^-$, and
\textbf{positive} if $U\subset (\NN)^+$.  Impose the following
transitive relation on multisets on $\bN^2$:
\begin{equation}
U\leq V\iff U_{(1)}-U_{(2)}\leq V_{(1)}-V_{(2)}.
\end{equation}

Let $\gi$ be the map on multisets on $\NN$ defined by
$\gi(\{(e_1,f_1),(e_2,f_2),\ldots\})=\{(f_1,e_1),(f_2,e_2),\ldots\}$.
Then $\iota$ is an involution, and it maps negative multisets on
$\NN$ to positive ones and visa-versa.  Thus $\iota$ is a
bijective pairing between the sets of negative and positive
multisets on $\NN$. Note also that $U\leq V\iff
\iota(V)\leq\iota(U)$.

In this paper, all sets and multisets other than $\bN$, $\NN$, and
multisets expressed explicitly as $A\disjunion (\bN\setminus C)$
where $A$ and $C$ are finite subsets of $\bN$, are assumed to be
finite.

\section{Semistandard Notched Bitableaux}\label{s.ss_notched_bitab}

A \textbf{notched bitableau} is a pair $(P, Q)$ of notched
tableaux of the same shape (i.e., the same number of rows and the
same number of boxes in each row).  The \textbf{degree} of $(P,Q)$
is the number of boxes in $P$ (or $Q$). A notched bitableau
$(P,Q)$ is said to be \textbf{row strict} if both $P$ and $Q$ are
row strict. A row strict notched bitableau $(P,Q)$ is said to be
\textbf{semistandard} if
\begin{equation}\label{e.s.ss_tabl_1}
P_1-Q_1\leq\cdots\leq P_r-Q_r.
\end{equation}
A row strict notched bitableau $(P,Q)$ is said to be
\textbf{negative} if $P_i\ltS Q_i$, $i=1,\ldots,r$,
\textbf{positive} if $P_i\gtS Q_i$, $i=1,\ldots,r$, and
\textbf{nonvanishing} if either
\begin{equation}\label{e.s.ss_tabl_2}
P_i\ltS Q_i\ \ \text{   or   }\ \ P_i\gtS Q_i,
\end{equation}
for each $i=1,\ldots,r$.

\begin{ex}
Consider the notched bitableau
$$(P,Q)=\left(\ \young(13,2357,2,346,45,6,467,79)\ \ \ ,
\ \ \ \young(59,4579,8,467,57,5,356,27)\right).$$ We have that
\begin{itemize}
\item[1.] $(P,Q)$ is row strict.

\item[2.] $P_1\disjunion Q_2=\{1,3,4,5,7,9\}\leq
\{2,3,5,5,7,9\}=P_2\disjunion Q_1$. Therefore, $P_1-Q_1\leq
P_2-Q_2$. Similarly, one checks that $P_i-Q_i\leq
P_{i+1}-Q_{i+1}$, $i=2,\ldots,7$.  Thus $(P,Q)$ is semistandard.

\item[3.] $P_i\ltS Q_i$, $i=1,\ldots,5$, and $P_i\gtS Q_i$,
$i=6,\ldots,8$. 
Thus $(P,Q)$ is nonvanishing.
\end{itemize}
\end{ex}

Let $(P,Q)$ be a semistandard notched bitableau. If, for subsets
$T$ and $W$ of $\bN^2$,
\begin{equation}\label{e.s.ss_tabl_on_TW}
T_{(1)}-T_{(2)}\leq P_1-Q_1 \quad\text{and}\quad P_r-Q_r\leq
W_{(1)}-W_{(2)},
\end{equation}
then we say that $(P,Q)$ is \textbf{bounded by} $\mathbf{T, W}$.
Note that (\ref{e.s.ss_tabl_on_TW}) combined with
(\ref{e.s.ss_tabl_1}) implies
\begin{equation*}
T_{(1)}-T_{(2)}\leq P_1-Q_1\leq\cdots\leq P_r-Q_r\leq
W_{(1)}-W_{(2)}.
\end{equation*}
Thus (\ref{e.s.ss_tabl_on_TW}) means that if $T_{(1)}$ and
$T_{(2)}$ are placed above the top rows of $P$ and $Q$
respectively, and $W_{(1)}$ and $W_{(2)}$ are placed below the
bottom rows of $P$ and $Q$ respectively, the the resulting
bitableau still satisfies (\ref{e.s.ss_tabl_1}).

Let $(P,Q)$ be a nonvanishing semistandard notched bitableau. Then
all rows $(P_i,Q_i)$ of $(P,Q)$ for which $P_i\ltS Q_i$ must lie
above all rows $(P_j,Q_j)$ of $(P,Q)$ such that $P_j\gtS Q_j$. Let
$0\leq i\leq r$ be maximal such that $P_i\ltS Q_i$. Then the top
$i$ rows of $P$ and $Q$ form a negative semistandard notched
bitableau and the bottom $r-i$ rows of $P$ and $Q$ form a positive
semistandard notched bitableau. These two bitableaux are called
respectively the \textbf{negative} and \textbf{positive parts} of
$(P,Q)$.

If $(P,Q)$ is a nonvanishing semistandard notched bitableau,
define $\gi(P,Q)$ to be the notched bitableau obtained by
reversing the order of the rows of $(Q,P)$. One checks that
$\gi(P,Q)$ is a nonvanishing semistandard notched bitableau. The
map $\gi$ is an involution, and it maps negative semistandard
notched bitableaux to positive ones and visa-versa. Thus $\gi$
gives a bijective pairing between the sets of negative and
positive semistandard notched bitableaux.

The definitions for semistandard notched tableaux and semistandard
notched bitableaux appear to be quite different. The following
Lemma gives a relationship between these two objects.

\begin{lem}\label{l.semist_bitab_implies_semist_on_b}
Let $(P,Q)$ be a negative semistandard notched bitableau, and let
$b$ be the minimum value of all entries of $Q$.  Then $P$ is
semistandard on $b$.
\end{lem}
\begin{proof}
Let $r$ be the number of rows of $P$. Suppose that, for some $i$,
$1\leq i\leq r-1$, $P_{i+1}$ has exactly $x$ entries which are
less than $b$, with $x>0$.  We must show that (i) $P_i$ has at
least $x$ entries which are less than $b$, and (ii)
$\{P_{i,1},\ldots,P_{i,x}\}\leq \{P_{i+1,1},\ldots,P_{i+1,x}\}$.

By (\ref{e.s.ss_tabl_1}),
\begin{equation}\label{e.ss_bitab_for_lem}
P_{i}\disjunion Q_{i+1}\leq P_{i+1}\disjunion Q_i.
\end{equation}
Therefore, since $P_{i+1}\disjunion Q_i$ has (exactly) $x$ entries
less than $b$, $P_{i}\disjunion Q_{i+1}$ must have at least $x$
entries less than $b$, which must all be in $P_i$, since $b$ is
the smallest entry of $Q$. This proves (i). The $x$ smallest
entries of the left hand side and right hand side of
(\ref{e.ss_bitab_for_lem}) are $\{P_{i,1},\ldots,P_{i,x}\}$ and
$\{P_{i+1,1},\ldots,P_{i+1,x}\}$ respectively.  Thus
(\ref{e.ss_bitab_for_lem}) implies (ii).
\end{proof}

\section{The Bounded RSK Correspondence}\label{s.grsk}

We next define the \textbf{bounded RSK correspondence}, $\rsk$, a
function which maps negative multisets on $\bN^2$ to negative
semistandard notched bitableaux. Let
$U=\{(a_1,b_1),\ldots,(a_t,b_t)\}$ be a negative multiset on
$\bN^2$, whose entries we assume are listed in
\textbf{lexicographic order}: (i) $b_1\geq\cdots\geq b_t$, and
(ii) if for any $i\in\{1,\ldots,t-1\}$, $b_i=b_{i+1}$, then
$a_i\geq a_{i+1}$. We inductively form a sequence of notched
bitableaux $(P^{(0)},Q^{(0)})$, $(P^{(1)},Q^{(1)}),$
$\ldots,(P^{(t)},Q^{(t)})$, such that $P^{(i)}$ is semistandard on
$b_i$, $i=1,\ldots,t$, as follows:
\begin{quote}
Let $(P^{(0)},Q^{(0)})=(\es,\es)$, and let $b_0=b_1$. Assume
inductively that we have formed $(P^{(i)}, Q^{(i)})$, such that
$P^{(i)}$ is semistandard on $b_i$, and thus on $b_{i+1}$, since
$b_{i+1}\leq b_i$. Define $P^{(i+1)}=P^{(i)}\la{b_{i+1}}a_{i+1}$.
Since bounded insertion preserves semistandardness on $b_{i+1}$,
$P^{(i+1)}$ is also semistandard on $b_{i+1}$. Let $j$ be the row
number of the new box of this bounded insertion. Define
$Q^{(i+1)}$ to be the notched tableau obtained by placing
$b_{i+1}$ on the {\it left} end of row $j$ of $Q^{(i)}$ (and
shifting all other entries of $Q^{(i)}$ to the right one box).
Clearly $P^{(i+1)}$ and $Q^{(i+1)}$ have the same shape.
\end{quote}
Then $\rsk(U)$ is defined to be $(P^{(t)},Q^{(t)})$. In the
process above, we write
$(P^{(i+1)},Q^{(i+1)})=(P^{(i)},Q^{(i)})\la{b_{i+1}}a_{i+1}$. In
terms of this notation,
$$\rsk(U)=((\es,\es)\la{b_1}a_1)\cdots\la{b_t}a_t.$$

\begin{ex}\label{ex.brsk}
Let $U=\{(7,8),(2,8),(6,7),(4,7),(1,7),(3,6),(2,4)\}$. Then
\[
\begin{array}{l@{\hspace{.9cm}}l}
P^{(0)}=\es
&Q^{(0)}=\es\\ \\
P^{(1)}=\es\la{8}7=\young(7)
&Q^{(1)}=\young(8)\\ \\
P^{(2)}=\young(7)\la{8}2=\young(2,7)
&Q^{(2)}=\young(8,8)\\ \\
P^{(3)}=\young(2,7)\la{7}6=\young(26,7)
&Q^{(3)}=\young(78,8)\\ \\
P^{(4)}=\young(26,7)\la{7}4=\young(24,67)
&Q^{(4)}=\young(78,78)\\ \\
P^{(5)}=\young(24,67)\la{7}1=\young(14,27,6)
&Q^{(5)}=\young(78,78,7)\\ \\
P^{(6)}=\young(14,27,6)\la{6}3=\young(13,247,6)
&Q^{(6)}=\young(78,678,7)\\ \\
P^{(7)}=\young(13,247,6)\la{4}2=\young(12,2347,6)
&Q^{(7)}=\young(78,4678,7)\\
\end{array}
\]
\vspace{.5em}

Therefore $\rsk(U)=\left(\, \young(12,2347,6)\ ,\
\young(78,4678,7)\right)$.
\end{ex}
\vspace{.5em}

\begin{lem}\label{l.ins_preserves_semist_bitab}
If $U$ is a negative multiset on $\bN^2$, then $\rsk(U)$ is a
negative semistandard notched bitableau.
\end{lem}
\begin{proof}
We use notation as in the definition of $\rsk$ above. That
$Q^{(i)}$ is row-strict for all $i$ follows from Lemma
\ref{l.bumping_route}(i) and the fact that the entries of $U$ are
listed in lexicographical order: if $b_{i+1}=b_i$ for some $i$,
then since $a_{i+1}\leq a_i$, the new box of the
$(i+1)^{\text{st}}$ insertion must be strictly below the new box
of $i^{\text{th}}$ insertion. That $P^{(i)}$ is row strict and
$(P^{(i)},Q^{(i)})$ is negative for all $i$ follows easily from
the definition of $\rsk$, using induction. It remains to prove
that $(P^{(i)},Q^{(i)})$ is semistandard for all $i$.

Let $P=P^{(i)}$, $Q=Q^{(i)}$, $P'=P^{(i+1)}$, $Q'=Q^{(i+1)}$,
$a=a_{i+1}$, $b=b_{i+1}$, and assume inductively that $(P,Q)$ is a
negative semistandard notched bitableau. Let $s$ be the number of
rows of $P'$ (and $Q'$). We show that $(P',Q')$ satisfies
(\ref{e.s.ss_tabl_1}), or equivalently, for any positive integer
$z$,
\begin{equation*} |(P'_j\disjunion
Q'_{j+1})^{<z}|\geq|(P'_{j+1}\disjunion Q'_{j})^{<z}|,\
j=1,\ldots,s-1.
\end{equation*}
By Lemma \ref{l.semist_bitab_implies_semist_on_b}, $P$ is
semistandard on $b$; hence so is $P'$.  Thus for $z\leq b$,
\begin{equation*}
|(P'_j\disjunion Q'_{j+1})^{<z}|=|(P'_j)^{< z}|\geq|(P'_{j+1})^{<
z}|=|(P'_{j+1}\disjunion Q'_{j})^{<z}|,\ j=1,\ldots,s-1.
\end{equation*}
Let $k$ be the row number of the new box (both in $P'$ and $Q'$)
of this bounded insertion. Since $(P,Q)$ is semistandard, for
$z>b$, $j\neq k-1$, and $j\neq k$,
\begin{equation*}
|(P'_j\disjunion Q'_{j+1})^{<z}|=|(P_j\disjunion
Q_{j+1})^{<z}|\geq|(P_{j+1}\disjunion Q_{j})^{<z}|
=|(P'_{j+1}\disjunion Q'_{j})^{<z}|.
\end{equation*}
For ($z>b$) and ($j= k-1$ or $j=k$),
\begin{equation*}
|(P'_j\disjunion Q'_{j+1})^{<z}|=|(P_j\disjunion
Q_{j+1})^{<z}|+1\geq|(P_{j+1}\disjunion Q_{j})^{<z}|+1
=|(P'_{j+1}\disjunion Q'_{j})^{<z}|.
\end{equation*}
\end{proof}

\begin{lem}\label{l.rsk_bijection__on_negatives}
The map $\rsk$ is a degree-preserving bijection from the set of
negative multisets on $\NN$ to the set of negative semistandard
notched bitableaux.
\end{lem}
\begin{proof}
That $\rsk$ is degree-preserving is obvious.

To show that $\rsk$ is a bijection, we define its inverse, which
we call the \textbf{reverse} of $\rsk$, or $\rrsk$.

Note that the bounded insertion used to form
$(P^{(i+1)},Q^{(i+1)})$ from $(P^{(i)},Q^{(i)})$,
$i=1,\ldots,t-1$, is reversible.  In other words, by knowing only
$(P^{(i+1)},Q^{(i+1)})$, we can retrieve $(P^{(i)},Q^{(i)})$,
$a_{i+1}$, and $b_{i+1}$. First, we obtain $b_{i+1}$; it is the
minimum entry of $Q^{(i+1)}$. Then, in the lowest row in which
$b_{i+1}$ appears in $Q^{(i+1)}$, select the greatest entry of
$P^{(i+1)}$ which is less than $b_{i+1}$. This entry was the new
box of the bounded insertion. If we begin reverse bounded
insertion with this entry, we retrieve $P^{(i)}$ and $a_{i+1}$.
Finally, $Q^{(i)}$ is retrieved from $Q^{(i+1)}$ by removing the
lowest occurrence of $b_{i+1}$ appearing in $Q^{(i+1)}$.  This
occurence must be on the left end of some row. All other entries
of that row should be moved one box to the left, thus eliminating
the empty box vacated by $b_{i+1}$.

It follows that we can reverse the entire sequence used to define
$\rsk$ by reversing each step in the sequence. If we generate
$(P^{(t)},Q^{(t)})$ via $\rsk$, we can retrieve $U$ using this
procedure. We will call the process of obtaining
$(P^{(i-1)},Q^{(i-1)})$, $a_i$, and $b_i$ from $(P^{(i)},Q^{(i)})$
described in the paragraph above a \textbf{reverse step} and
denote it by $(P^{(i-1)},Q^{(i-1)})=(P^{(i)},Q^{(i)})\ra{b_i}a_i$.
 We will call the process of applying all the reverse steps sequentially
to retrieve $U$ from $(P^{(t)},Q^{(t)})$ the \textbf{reverse of}
$\rsk$, or $\rrsk$.  For example, if one applies $\rrsk$ to the
negative semistandard notched bitableau appearing on the bottom
line of Example \ref{ex.brsk}, one obtains the negative multiset
$U$ from that example.


If $(P^{(t)},Q^{(t)})$ is an arbitrary semistandard notched
bitableau (which we do not assume to be $\rsk(U)$, for some $U$),
then we can still apply a sequence of reverse steps to
$(P^{(t)},Q^{(t)})$, to sequentially obtain $(P^{(i)},Q^{(i)})$,
$a_i$, $b_i$, $i=t,\ldots,1$. For this process to be well-defined,
however, it must first be checked that the successive
$(P^{(i)},Q^{(i)})$ are negative semistandard notched bitableaux.
It suffices to prove a statement very similar to that proved in
Lemma \ref{l.ins_preserves_semist_bitab}: if $(P,Q)$ is a negative
semistandard notched bitableau, then $(P',Q'):=(P,Q)\ra{b}a$ is a
negative semistandard notched bitableau, $a<b$ are positive
integers, and $b$ is less than or equal to all entries of $Q$.
That $a<b$ are positive integers and $b$ is less than or equal to
all entries of $Q$ follow immediately from the definition of a
reverse step. That $(P',Q')$ is a negative semistandard notched
bitableau follows in much the same manner as the proof of Lemma
\ref{l.ins_preserves_semist_bitab}; we omit the details.

It remains to show that the elements
$\{(a_1,b_1),\ldots,(a_t,b_t)\}$ of the negative multiset on
$\bN^2$ produced by applying this sequence of reverse steps to the
arbitrary semistandard notched bitableau $(P^{(t)},Q^{(t)})$ are
listed in lexicographic order. That $b_{i}\geq b_{i+1}$ follows
from the definition of $\rrsk$: $b_{i+1}$ is the minimim entry of
$Q^{(i+1)}$, which also has $b_i$ as an entry. If $b_i=b_{i+1}$,
then $a_i\geq a_{i+1}$ is a consequence of Lemma
\ref{l.bumping_route}(i) and (ii).

At each step, $\rsk$ and the reverse of $\rrsk$ are inverse to
eachother. Thus they are inverse maps.

\end{proof}

\begin{figure}[h!]
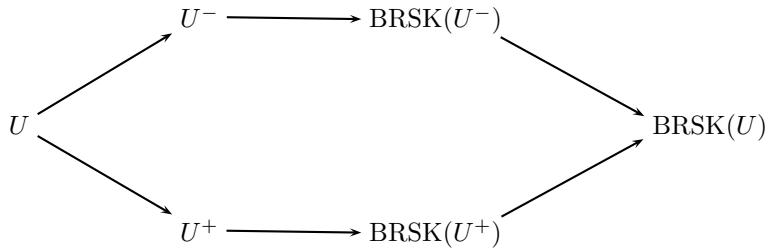

\begin{center}
$
\setlength{\arraycolsep}{1cm}
\begin{array}{llll}
&\rnode{b}{U^-} &\rnode{c}{\rsk(U^-)}&\\[1cm]
\rnode{a}{U}&&& \rnode{d}{\rsk(U)}\\[1cm]
&\rnode{B}{U^+} & \rnode{C}{\rsk(U^+)}&
\end{array}
\ncline[nodesep=3pt]{->}{a}{b}
\ncline[nodesep=3pt]{->}{b}{c}
\ncline[nodesep=4pt,offsetA=3pt,offsetB=-5pt]{->}{c}{d}
\ncline[nodesep=3pt]{->}{a}{B}
\ncline[nodesep=3pt]{->}{B}{C}
\ncline[nodesep=4pt,offsetA=-3pt,offsetB=5pt]{->}{C}{d}
$
\end{center}
\caption{\label{f.the_map_ersk}The map $\rsk$}
\end{figure}

The map $\rsk$ can be extended to all nonvanishing multisets on
$\bN^2$. If $U$ is a positive multisets on $\bN^2$, then define
$\rsk(U)$ to be $\gi(\rsk(\gi(U)))$. If $U$ is a nonvanishing
multisets on $\bN^2$, with negative and positive parts $U^-$ and
$U^+$, then define $\rsk(U)$ to be the semistandard notched
bitableau whose negative and positive parts are $\rsk(U^-)$ and
$\rsk(U^+)$ (see Figure \ref{f.the_map_ersk}). \vspace{10pt}
As a consequence of Lemma \ref{l.rsk_bijection__on_negatives}, we
obtain
\begin{prop}\label{p.rsk_bijection}
The map $\rsk$ is a degree-preserving bijection from the set of
nonvanishing (resp. negative, positive) multisets on $\bN^2$ to
the set of nonvanishing (resp. negative, positive) semistandard
notched bitableaux.
\end{prop}

The ordinary RSK correspondence is a degree-preserving bijection
from the set of multisets on $\bN^2$ to the set of semistandard
bitableaux. The process used to define the bijection is similar to
the process described above to define the BRSK. There are two
essential differences between the two processes. First, in the
ordinary RSK, the multiset is not first seperated into its
positive and negative parts. Indeed, the ordinary RSK is oblivious
as to whether elements of the multiset are positive or negative.
Secondly, in the ordinary RSK, ordinary Schensted insertion is
used rather than bounded insertion. See \cite{Fu} or \cite{Sa} for
more details on the ordinary RSK.

\section{Restricting the Bounded RSK Correspondence}\label{s.grsk_props}

Thus far, there has been no reference to $\ga$, $\gb$, or $\gc$ in
our definition or discussion of the bounded RSK. In fact, each of
$\ga$, $\gb$, and $\gc$ is used to impose restrictions on the
domain and codomain of the bounded RSK. It is the bounded RSK,
with domain and codomain restricted according to $\ga$, $\gb$, and
$\gc$, which is used in Section \ref{s.grobner} to give
geometrical information about $Y\abc$.

In this section, we first show how $\gb$ restricts the domain and
codomain of the bounded RSK.  We then show how two subsets $T$ and
$W$ of $\NN$, $T$ negative and $W$ positive, restrict the domain
and codomain of the bounded RSK. In Section \ref{s.grobner}, these
two subsets will be replaced by $T_\ga$ and $W_\gc$, subsets of
$\NN$ determined by $\ga$ and $\gc$ respectively.

\subsection*{Restricting by $\gb$}

Let $\gb\in I_{d,n}$. We say that a notched bitableau $(P,Q)$ is
\textbf{on} $\gB$ if all entries of $P$ are in $\ol{\gb}$ and all
entries of $Q$ are in $\gb$. It is clear from the construction of
$\rsk$ that if $U$ is a nonvanishing multiset on $\gB$, then
$\rsk(U)$ is a nonvanishing semistandard notched bitableau on
$\gB$, and visa-versa. Thus, as a consequence of Corollary
\ref{p.rsk_bijection}, we obtain

\begin{cor}\label{c.beta_rsk_bijection}
The map $\rsk$ restricts to a degree-preserving bijection from the
set of nonvanishing (resp. negative, positive) multisets on $\BBZ$
to the set of nonvanishing (resp. negative, positive) notched
bitableaux on $\BBZ$.
\end{cor}
\noindent We remark that if $\gb$ is the largest or smallest
element of $I_{d,n}$ ($\{n-d+1,\ldots,n\}$ or $\{1,\ldots,d\}$
respectively), then the bounded RSK restricted to $\gB$ is the
same algorithm as the ordinary RSK restricted to $\gB$.

\subsection*{Restricting by $T$ and $W$}

A \textbf{chain} in $\NN$ is a subset
$C=\{(e_1,f_1),\ldots,(e_m,f_m)\}$ of $\NN$ such that
$e_1<\cdots<e_m$ and $f_1>\cdots>f_m$.
Let $T$ and $W$ be negative and positive subsets of $\NN$
respectively. A nonempty multiset $U$ on $\NN$ is said to be
\textbf{bounded by} $\mathbf{T, W}$ if for every chain $C$ which
is contained in the underlying set of $U$,
\begin{equation}\label{e.g.admissible_def}
T\leq C\leq W
\end{equation}
(where we use the order on multisets on $\bN^2$ defined in Section
\ref{s.multisets}).  We note that, in general, this condition
neither implies nor is implied by the condition $T\leq U\leq W$.
For special cases, a geometric interpretation in terms of a {\it
chain order} for $U$ being bounded by $T,W$ appears in Section
\ref{s.twisted} (this interpretation is not necessary for our
discussion here).

With this definition, the bounded RSK correspondence is a bounded
function, in the sense that it maps bounded sets to bounded sets.
More precisely, we have the following Lemma, whose proof appears
in Section \ref{s.proofs}.

\begin{lem}\label{l.g.dominance}
If a nonvanishing multiset $U$ on $\NN$ is bounded by $T,W$, then
$\rsk(U)$ is bounded by $T,W$.
\end{lem}

Let $T$ and $W$ be negative and positive subsets of $\gB$,
respectively. Combining Corollary \ref{c.beta_rsk_bijection} and
Lemma \ref{l.g.dominance}, we obtain
\begin{cor}\label{c.g.beta_dominance_2}
For any positive integer $m$, the number of degree $m$
nonvanishing multisets on $\BBZ$ bounded by $T, W$ is less than or
equal to the number of degree $m$ nonvanishing semistandard
notched bitableaux on $\BBZ$ bounded by $T, W$.
\end{cor}

\noindent We mention that the converse of Lemma
\ref{l.g.dominance} is not a priori true, i.e., the reverse $\rsk$
is not a priori bounded.  Otherwise, we could state here that the
two numbers in Corollary \ref{c.g.beta_dominance_2} are equal.  In
fact, the reverse $\rsk$ is indeed bounded, but since we do not
need this result for our purposes we omit the proof.

\section{A Gr\"obner Basis}\label{s.grobner}

We call $f=f_{\gt_1,\gb}\cdots f_{\gt_r,\gb}\in K[\cO_\gb]$ a
\textbf{standard monomial} if
\begin{equation}\label{e.r.stand_mon_1}
\gt_1\leqT \cdots\leqT \gt_r
\end{equation}
and for each $i\in\{1,\ldots,r\}$, either
\begin{equation}\label{e.r.stand_mon_2}
\gt_i\ltT \gb\ \ \hbox{ or }\ \ \gt_i\gtT \gb.
\end{equation}
If in addition, for $\ga,\gc\in I_{d,n}$,
\begin{equation}\label{e.r.stand_mon_on_ab} \ga\leq \gt_1\quad\text{and} \quad
\gt_r\leq\gc,
\end{equation}
then we say that $f$ is \textbf{standard on} $Y\abc$.

We remark that, in general, a standard monomial is not a monomial
in the affine coordinates of $\cO_\gb$ (the $x_{ij}$'s); rather,
it is a polynomial. It is only a monomial in the $f_{\gt,\gb}$'s.
Recall the following result (see \cite{Kr-La3}).
\begin{thm}\label{t.r.stand_mons_form_basis}
The standard monomials on $Y_{\ga,\gb}^\gc$ form a basis for
$K[Y_{\ga,\gb}^\gc]$.
\end{thm}

We wish to give a different indexing set for the standard
monomials on $Y_{\ga,\gb}^\gc$. Let $I_\gb$ be the pairs
$(\pr,\qr)$ such that $\pr\subset\ol{\gb}$, $\qr\subset\gb$, and
$|\pr|=|\qr|$.  Defining $\pr-\qr:=\pr\disjunion(\gb\setminus
\qr)$ (see Section \ref{s.multisets}), we have the following fact,
which is easily verified:
\begin{equation*}
\hbox{The map }(\pr,\qr)\mapsto \pr-\qr\hbox{ is a bijection from
}I_\gb\hbox{ to }I_{d,n},
\end{equation*}
(Indeed, the inverse map is given by $\gt\mapsto
(\gt\setminus\gb,\gb\setminus\gt)$). Thus, for instance, we may
write a Plucker coordinate in $\bC[\cO_\gb]$ as $f_{R-S}$ instead
of $f_{\gt,\gb}$. (In fact, $f_{R-S}$ is just plus or minus the
minor of $\bC[\cO_\gb]$ with row set $R$ and column set $S$; in
Example \ref{ex.big-cell},
$f_{\{1,4\}-\{2,7\}}=f_{\{1,4,5\},\gb}$.)

Note that under this bijection, $(\es,\es)$ maps to $\gb$. Let
$(\pr_\ga,\qr_\ga)$ and $(\pr_\gc,\qr_\gc)$ be the preimages of
$\ga$ and $\gc$ respectively.  Define $T_\ga$ and $W_\gc$ to be
any subsets of $\BBZ$ such that $(T_\ga)_{(1)}=\pr_\ga$,
$(T_\ga)_{(2)}=\qr_\ga$, $(W_\gc)_{(1)}=\pr_\gc$,
$(W_\gc)_{(2)}=\qr_\gc$.

Given any notched bitableau $(P,Q)$ which is on $\BBZ$, we can
form the monomial $f=f_{P_1-Q_1}\cdots f_{P_r-Q_r}\in K[\cO_\gb]$,
where $P_1,\ldots,P_r$ are the rows of $P$ and $Q_1,\ldots,Q_r$
the rows of $Q$. Conversely, given any monomial of the form
$f=f_{P_1-Q_1}\cdots f_{P_r-Q_r}\in K[\cO_\gb]$, where $P_i$ and
$Q_i$ are subsets of $\ol{\gb}$ and $\gb$ respectively of the same
cardinality, $i=1,\ldots,r$, we can form the notched bitableau
$(P,Q)$ on $\BBZ$ whose rows are $(P_1,Q_1),\ldots,(P_r,Q_r)$.  We
have
\begin{align*}
(P,Q)&\text{ is semistandard, nonvanishing, and bounded by }
T_\ga,
W_\gc\\
&\iff  (P,Q) \text{ satisfies (\ref{e.s.ss_tabl_1}),
(\ref{e.s.ss_tabl_2}), and (\ref{e.s.ss_tabl_on_TW})}\\
&\iff (P_1,Q_1),\ldots,(P_r,Q_r) \text{ satisfy
(\ref{e.s.ss_tabl_1}),
(\ref{e.s.ss_tabl_2}), and (\ref{e.s.ss_tabl_on_TW})}\\
&\iff P_1-Q_1,\ldots,P_r-Q_r \text{ satisfy
(\ref{e.r.stand_mon_1}),
(\ref{e.r.stand_mon_2}), and (\ref{e.r.stand_mon_on_ab})}\\
&\iff f_{P_1-Q_1}\cdots f_{P_r-Q_r} \text{ is standard on } Y\abc.
\end{align*}
When we write above that $(P_1,Q_1),\ldots,(P_r,Q_r)$ satisfy
(\ref{e.s.ss_tabl_1}) and (\ref{e.s.ss_tabl_on_TW}), we use
(\ref{e.m.set_subtraction}) to describe the order.  When we write
that $P_1-Q_1,\ldots,P_r-Q_r$ satisfy (\ref{e.r.stand_mon_1}) and
(\ref{e.r.stand_mon_on_ab}), we use the termwise order on the
$P_{i}\disjunion(\gb\setminus Q_{i})$. The equivalence of these
orders is discussed in Section \ref{s.multisets}. To see that
$(P_i,Q_i)$ satisfying (\ref{e.s.ss_tabl_2}) is equivalent to
$P_i-Q_i$ satisfying (\ref{e.r.stand_mon_2}), note that since
$\ol{\gb}\cap\gb=\es$, $P_i\cap Q_i=\es$, and therefore
$P_i\ltS Q_i\iff P_i-Q_i=P_i\disjunion (\gb\setminus Q_i)<\gb$.

This proves the following lemma.

\begin{lem}\label{l.r.semist_index}
The degree $m$ nonvanishing semistandard notched bitableaux on
$\BBZ$ bounded by $T_\ga, W_\gc$ form an indexing set for the
degree $m$ standard monomials on $Y_{\ga,\gb}^\gc$.
\end{lem}

As discussed in Section \ref{s.results}, the affine coordinates of
$K[\cO_\gb]$ are indexed by $\BBZ$.  Thus monomials in the affine
coordinates of $K[\cO_\gb]$ are naturally indexed by multisets on
$\BBZ$: the monomial $x_{i_1j_1}\ldots x_{i_tj_t}\in K[\cO_\gb]$
is indexed by the multiset
$\{(i_1,j_1),\ldots,(i_t,j_t)\}\subset\BBZ$. Letting
$U=\{(i_1,j_1),\ldots,(i_t,j_t)\}\subset\BBZ$, we shall denote
$x_{i_1j_1}\ldots x_{i_tj_t}$ by $x_U$. Note that $x_U$ is
square-free if and only if the multiset $U$ is in fact a set. We
define a monomial order on $K[\cO_\gb]$ as follows: first declare
$x_{ij}<x_{i'j'}$ if ($i<i'$) or ($i=i'$ and $j>j'$), then impose
the lexicographic order on monomials. For our purposes, the
critical feature of this monomial order is that the initial term
of any minor $f_{R-S}$ is the Southwest-Northeast monomial of
$f_{R-S}$ (when $f_{R-S}$ is written naturally as plus or minus a
determinant all of whose entries are $x_{ij}$'s). In other words,
$\init f_{R-S}=x_C$, where $C$ is the chain in $\gB$ with
$C_{(1)}=R$, $C_{(2)}=S$ (in Example \ref{ex.big-cell}, $\init
f_{\{1,4\}-\{2,7\}}=x_{17}x_{42}=x_{\{(1,7),(4,2)\}}$). Any other
monomial order with this property would also suit our purposes.

\begin{lem}\label{l.r.admiss_index}
The degree $m$ multisets on $\BBZ$ bounded by $T_\ga, W_\gc$ form
an indexing set for the degree $m$ monomials of
$K[\cO_\gb]\setminus \init G\abc$.
\end{lem}
\begin{proof}
\begin{align*}
\init G\abc&= \langle\init f_{\gt,\gb}\mid \ga\not\leq\gt\hbox{ or
} \gt\not\leq\gc\rangle\\
&=\langle\init f_{\pr-\qr}\mid \pr_\ga-\qr_\ga\not\leq
\pr-\qr\hbox{ or }
\pr-\qr\not\leq \pr_\gc-\qr_\gc\rangle\\
&=\langle x_C\mid C\hbox{ a chain}, \pr_\ga-\qr_\ga\not\leq
C_{(1)}-C_{(2)}\hbox{ or } C_{(1)}-C_{(2)}\not\leq
\pr_\gc-\qr_\gc\rangle\\
&=\langle x_C\mid C\hbox{ a chain}, T_\ga\not\leq C \hbox{ or }
C\not\leq W_\gc\rangle.
\end{align*}
Therefore,\\

\noindent $x_U$ is a monomial in $K[\cO_\gb]\setminus \init
G_{\ga,\gb}^\gc$
\begin{itemize}
\item[$\iff$] $x_U$ is not divisible by any $x_C$, $C$ a chain
such that $T_\ga\not\leq C$ or $C\not\leq W_\gc$

\item[$\iff$] $U$ contains no chains $C$ such that $T_\ga\not\leq
C$ or $C\not\leq W_\gc$

\item[$\iff$] $T_\ga\leq C\leq W_\gc$, for any chain $C$ in $U$

\item[$\iff$] $U$ is bounded by $T_\ga,W_\gc$.
\end{itemize}
\end{proof}
\noindent We are now ready to prove the main result of the paper.

\begin{proof}[Proof of Proposition \ref{p.i.grobner}.] We wish to show that $\init G_{\ga,\gb}^\gc= \init
\langle G\abc\rangle$.  Since $G_{\ga,\gb}^\gc\subset \langle
G\abc\rangle$,
$\init G_{\ga,\gb}^\gc\subset \init \langle G\abc\rangle$. For any $m\geq 1$,\\

\noindent\# of degree $m$ monomials in $K[\cO_\gb]\setminus \init
G_{\ga,\gb}^\gc$
\begin{itemize}
\item[$\stackrel{a}{=}$] \# of degree $m$ multisets on $\BBZ$
bounded by $T_\ga,W_\gc$

\item[$\stackrel{b}{\leq}$] \# of degree $m$ semistandard notched
bitableaux on $\BBZ$ bounded by $T_\ga,W_\gc$

\item[$\stackrel{c}{=}$] \# of degree $m$ standard monomials on
$Y_{\ga,\gb}^\gc$

\item[$\stackrel{d}{=}$] \# of degree $m$ monomials in
$K[\cO_\gb]\setminus \init \langle G\abc\rangle$,
\end{itemize}

\noindent where $a$ follows from Lemma \ref{l.r.admiss_index}, $b$
from Corollary \ref{c.g.beta_dominance_2}, $c$ from Lemma
\ref{l.r.semist_index}, and $d$ from the fact that standard
monomials on $Y_{\ga,\gb}^\gc$ and the monomials in
$K[\cO_\gb]\setminus \init \langle G\abc\rangle$ both induce
homogeneous bases for $K[\cO_\gb]/ \langle G\abc\rangle$. Thus
$\init G_{\ga,\gb}^\gc\supset \init \langle G\abc\rangle$.

We point out that, as a consequence of this proof, inequality b is
actually an equality.
\end{proof}

\begin{ex}
Let $n=6$, $d=3$, $\gb=\{3,5,6\}$. Then
\begin{align*}
\cO_\gb&=\left\{ \left(
\begin{matrix}
x_{13}&x_{15}&x_{16}\\
x_{23}&x_{25}&x_{26}\\
1&0&0\\
x_{43}&x_{45}&x_{46}\\
0&1&0\\
0&0&1
\end{matrix}
\right), x_{ij}\in K \right\}
\end{align*}
Let $\ga=\{1,2,4\}$, $\gc=\{4,5,6\}$. We list all objects
identified with the monomial
$x_{26}\,x_{45}^2\,x_{15}\,x_{13}\,x_{43}\in K[\cO_\gb]\setminus
\init G_{\ga,\gb}^\gc$ in the preceding discussion: (a) monomial
in $K[\cO_\gb]\setminus \init G_{\ga,\gb}^\gc$, (b) multiset on
$\BBZ$ bounded by $T_\ga,W_\gc$, (c) semistandard notched
bitableau on $\BBZ$ bounded by $T_\ga,W_\gc$, (d) standard
monomial on $Y_{\ga,\gb}^\gc$, and (e) standard monomial on
$Y_{\ga,\gb}^\gc$ (different indices).
\begin{itemize}
\item[(a)] $x_{26}\,x_{45}^2\,x_{15}\,x_{13}\,x_{43}$
\item[(b)] $\{(2,6),(4,5),(4,5),(1,5),(1,3),(4,3)\}$
%
\item[(c)] $\left(\ \young(14,1,24,4)\ , \ \young(56,5,35,3)\
\right)$
\item[(d)]
$f_{\{1,4\}-\{5,6\}}\,f_{\{1\}-\{5\}}\,f_{\{2,4\}-\{3,5\}}\,f_{\{4\}-\{3\}}$
\item[(e)]
$f_{\{1,3,4\},\gb}\,f_{\{1,3,6\},\gb}\,f_{\{2,4,6\},\gb}\,f_{\{4,5,6\},\gb}$.
\end{itemize}
Note that the semistandard notched tableau in (c) is obtained from
the multiset in (b) by applying the BRSK.
\end{ex}

\noindent Consider the following general Lemma on Gr\"obner Bases
(see \cite{Ei}).
\begin{lem}\label{l.r.square_free_grobner}
Let $R=K[x_1,\ldots,x_m]$ be a polynomial ring, let $I\subset R$
be a homogeneous ideal, and let $G=\{g_1,\ldots,g_k\}$ be a
Gr\"obner basis for $I$, such that $\init(g_i)$ is square-free,
$i=1,\ldots,k$.  Let $M$ be the maximum degree of a square-free
monomial in $R\setminus\init(G)$. Then $\dim(R/I)=M$, and
$\deg(R/I)$ is the number of square-free monomials of degree $M$
in $R\setminus \init(G)$.
\end{lem}

\noindent Since the initial term of $f_{R-S}$ is square-free for
any $(\pr,\qr)\in I_\gb$, Lemma \ref{l.r.square_free_grobner} may
be applied to our situation in order to obtain Corollary
\ref{c.i.grobner_extras}. Indeed, by Lemma
\ref{l.r.square_free_grobner}, $\deg(Y\abc)$ is the number of
square-free monomials of degree $M$ in $K[\cO_\gb]\setminus \init
G\abc$, where $M=\dim(Y\abc)=l(\gc)-l(\ga)$. We remark that, in
Lemma \ref{l.r.square_free_grobner}, we use the convention that 1
is the only square-free monomial of degree zero.

\section{Twisted Chains and Multiplicities}\label{s.twisted}

The goal of this section and the next one is to establish
Proposition \ref{p.mult_chain}, which also appears in
\cite{Ko-Ra}, \cite{Kra1}, \cite{Kra2}, and \cite{Kr-La}.
Proposition \ref{p.mult_chain} gives a combinatorial formula for
multiplicities which involves counting families of nonintersecting
lattice paths. Proposition \ref{p.mult_chain} is essentially a
reformulation of Corollary \ref{c.i.grobner_extras} in more
combinatorial language. We establish Proposition
\ref{p.mult_chain} in two steps. In this section, we show that
Corollary \ref{c.i.grobner_extras} implies Lemma
\ref{l.chain_order}; in Section \ref{s.path_families_mults}, we
show that Lemma \ref{l.chain_order} implies Proposition
\ref{p.mult_chain}.

We begin by introducing {\it twisted chains} and {\it chain
boundedness}, notions which allow us to place earlier results on
combinatorial footing.  We define the following partial orders on
the negative elements of $\NN$: if $(e,f),(g,h)\in\NN$, both
negative, then
\begin{align*}
&(e,f)\ltC (g,h)\ \ \hbox{ if }\ \ f< h\hbox{ and }e> g,\\
&(e,f)\leqCC (g,h)\ \ \hbox{ if }\ \ f\leq h\hbox{ and }e\geq g.
\end{align*}
Note that $\leqCC$ is a finer order than $\ltC$. If
$(c,d),(e,f)\in(\NN)^-$, then define $$(c,d)\meet
(e,f)=(\max(c,e),\min(d,f))\in\NN.$$  If
$T=\{(e_1,e_2),\ldots,(e_m,e_{m+1})\}$ is a subset of $\NN$, then
we say that $T$ is \textbf{completely disjointed} if $e_i\neq e_j$
when $i\neq j$. A \textbf{negative twisted chain} is a completely
disjointed negative subset of $\NN$ such that for any $u,v\in T$,
$u\neq v$, either $u\ltC v$, $v\ltC u$, or $u\meet
v\not\in(\NN)^-$.


\begin{ex}
A negative chain in $\NN$, defined in Section \ref{s.grsk_props},
can alternatively be described as a negative subset
$\{u_1,\ldots,u_m\}$ of $\NN$ such that $u_1\ltC\cdots\ltC u_m$. A
negative chain is a negative twisted chain.
\end{ex}

Let $T=\{(e_1,f_1),$ $\ldots,(e_m,f_m)\}$ be a completely
disjointed negative subset of $\NN$ such that $f_1<\cdots<f_t$.
For $\gs\in S_m$, the permutation group on $m$ elements, we define
$\gs(T)=\{(e_{\gs(1)},f_1),\ldots,(e_{\gs(m)},f_m)\}$. Let
$\cT=\{\gs(T)\mid \gs\in S_n, \gs(T)\,\text{negative}\}$. Impose
the following total order on $\cT$: if
$R=\{(a_1,f_1),\ldots,(a_t,f_t)\},$ $
S=\{(b_1,f_1),\ldots,(b_t,f_t)\}$ $\in\cT$, then
$R\stackrel{\hbox{\tiny lex}}{<} S$ if, for the smallest $i$ for
which $a_i\neq b_i$, $a_i>b_i$. Since $\stackrel{\hbox{\tiny
lex}}{<}$ is a total order, $\cT$ has a unique minimal element,
which we denote by $\wt{T}$.

\begin{lem}\label{l.wtt_twisted_chain}
$\wt{T}$ is a negative twisted chain.
\end{lem}
\begin{proof}
Suppose that $\wt{T}=\{(c_1,f_1),\ldots,(c_m,f_m)\}$ is not a
negative twisted chain.  Then there exists $i<j$ such that
$(c_i,f_i)\not\ltC (c_j,f_j)$, $(c_j,f_j)\not\ltC (c_i,f_i)$, and
$(c_i,f_i)\meet (c_j,f_j)\in(\NN)^-$. This implies that
$c_i<c_j<f_i<f_j$. Letting $\gs_{i,j}$ be the transposition which
switches $i$ and $j$, we have that $\gs_{i,j}(\wt{T})$
$=\{(c_1,f_1),\ldots,$
$(c_{j},f_i),\ldots,(c_i,f_{j}),\ldots,(c_m,f_m)\}$. Since
$c_{j}<f_i$ and $c_i<f_{j}$, $\gs_{i,j}(\wt{T})$ is negative, and
hence $\gs_{i,j}(\wt{T})\in\cT$. The fact that $c_i<c_{j}$ implies
that $\gs_{i,j}(\wt{T})\stackrel{\hbox{\tiny lex}}{<} \wt{T}$,
which contradicts the minimality of $\wt{T}$.
\end{proof}

\eject

\begin{ex} The set $\NN$ is plotted in both (a) and (b) below.  The large squares cover
the points of $\NN\setminus (\NN)^{\neq 0}$, and thus separate the
points of $(\NN)^-$ from those of $(\NN)^+$. In (a), the
$\times$'s form a completely disjointed negative subset $T$ of
$\NN$. In (b), the $\times$'s form $\wt{T}$. Note that $\wt{T}$ is
a negative twisted chain, as required by Lemma
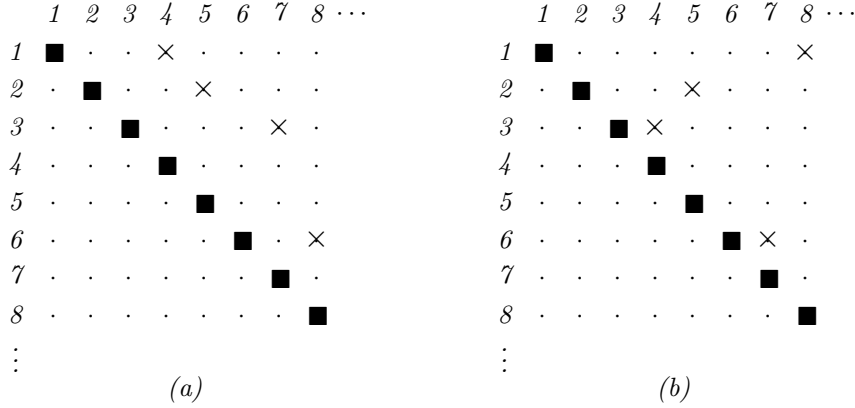
\ref{l.wtt_twisted_chain}.

\begin{center}
\setlength{\unitlength}{.5mm}
$\begin{array}{l@{\hspace{1.5cm}}l}
\begin{picture}(100,110)(0,0)
\matrixput(10,20)(10,0){8}(0,10){8}{\plotsmcirchar}
\multiputlist(0,10)(0,10){$\vdots$,8,7,6,5,4,3,2,1}
\multiputlist(10,100)(10,0){1,2,3,4,5,6,7,8,$\cdots$}
\multiput(10,90)(10,-10){8}{\plotsmbcirchar}

\put(40,90){\plotstarchar} \put(50,80){\plotstarchar}
\put(70,70){\plotstarchar} \put(80,40){\plotstarchar}

\put(45,0){\makebox(0,0){\textrm{(a)}}}

\end{picture}
&

\begin{picture}(100,110)(0,0)
\matrixput(10,20)(10,0){8}(0,10){8}{\plotsmcirchar}
\multiputlist(0,10)(0,10){$\vdots$,8,7,6,5,4,3,2,1}
\multiputlist(10,100)(10,0){1,2,3,4,5,6,7,8,$\cdots$}
\multiput(10,90)(10,-10){8}{\plotsmbcirchar}

\put(40,70){\plotstarchar} \put(50,80){\plotstarchar}
\put(70,40){\plotstarchar} \put(80,90){\plotstarchar}

\put(45,0){\makebox(0,0){\textrm{(b)}}}

\end{picture}

\end{array}
$
\end{center}
\end{ex}

\vspace{1em}

If $T$ is a positive subset of $\NN$, then we say that $T$ is a
\textbf{positive twisted chain} if $\gi(T)$ is a negative twisted
chain. A \textbf{twisted chain} is a subset of $\NN$ which is
either a positive or a negative twisted chain.

For $R$ a negative subset of $\NN$ and $x\in \NN$ negative, define
$\mathbf{\text{\textbf{depth}}_R(x)}$ to be the maximum $r$ such
that there exists a chain $u_1\ltC\cdots\ltC u_r$ in $R$ with
$u_r\leqCC x$. We extend $\leqCC$ to a transitive relation on
subsets of $\bN^2$ as follows. If $R$, $S$ are negative subsets of
$\bN^2$, then $R\leqCC S$ (or $S\geqCC R$) if
$\depth_R(x)\geq\depth_S(x)$ for every negative $x\in\bN^2$.  Note
that this is equivalent to $\depth_R(x)\geq\depth_S(x)$ for every
$x\in S$. 
If $R$, $S$ are positive subsets of $\bN^2$, then we say that
$S\geqCC R$ if $\gi(S)\leqCC \gi(R)$. If $R$ is a negative subset
of $\bN^2$ and $S$ is a positive subset, then we say that $R\leqCC
S$.

Recall the relation $\leq$ on multisets on $\bN^2$ defined in
Section \ref{s.multisets}. The following Lemma, whose proof
appears in Section \ref{s.proofs}, provides the key step in the
proof of Lemma \ref{l.chain_order}.

\begin{lem}\label{l.t.two_orders}
Let $R$ and $S$ be twisted chains. Then $R\leqCC S\iff R\leq S$.
\end{lem}

Let $R$ and $S$ be negative and positive twisted chains
respectively. We say that a multiset $U$ on $\NN$ is
\textbf{chain-bounded by} $\mathbf{R,S}$ if $R\leqCC U^-$ and
$U^+\leqCC S$, or equivalently, if for every chain $C$ in $U$,
\begin{equation*}
R\leqCC C^-\ \ \text{ and }\ \ C^+\leqCC S.
\end{equation*} In (\ref{e.g.admissible_def}), one can replace
$R\leq C\leq S$ by $R\leq C^-$ and $C^+\leq S$. Thus, by Lemma
\ref{l.t.two_orders}, $U$ is chain-bounded by $R,S$ if and only if
$U$ is bounded by $R,S$.

For the remainder of this section and the next one, we will be
interested in twisted chains which are contained in $\gB$, a
subset of $\NN$. 
Example \ref{e.chain_order} illustrates two negative twisted
chains in $\gB$.

\begin{ex}\label{e.chain_order}
The set $\gB\subset\NN$, for $d=8$, $n=17$,
$\gb=\{2,7,8,9,12,13,16,17\}$, is plotted below.  The dotted line
separates the negative from positive elements of $\gB$. The
$\times$'s form a negative twisted chain $R$ in $\gB$; the
$\bullet$'s form a negative twisted chain $S$ in $\gB$; and
$R\leqCC S$.

\begin{center}
\setlength{\unitlength}{.5mm}
\begin{picture}(90,110)(0,0)
\matrixput(10,10)(10,0){8}(0,10){9}{\plotsmcirchar}
\multiputlist(0,10)(0,10){15,14,11,10,6,5,4,3,1}
\multiputlist(10,100)(10,0){2,7,8,9,12,13,16,17}

\dottedline{.7}(10,85)(15,85)(15,45)(45,45)(45,25)(65,25)(65,10)

\put(80,90){\plotstarchar} \put(70,20){\plotstarchar}
\put(60,80){\plotstarchar} 
\put(20,70){\plotstarchar}

\put(80,20){\plotcirchar} \put(70,10){\plotcirchar}
\put(60,40){\plotcirchar} \put(50,30){\plotcirchar}
\put(40,80){\plotcirchar} \put(30,60){\plotcirchar}
\put(20,50){\plotcirchar} \put(10,90){\plotcirchar}
\end{picture}
\end{center}

\end{ex}

\begin{lem}\label{l.chain_order}
$\Mult_{e_\gb}X_\ga^\gc$ is the number of subsets $U$ of $\BBZ$
which are of maximal degree  among those which are chain-bounded
by $\wt{T_\ga}, \wt{W_\gc}$.
\end{lem}
\begin{proof}
Recall that if $U$ is a multiset on $\BBZ$, then the monomial
$x_U$ is square-free if and only if $U$ is a subset of $\BBZ$,
i.e., each of its elements has degree 1.  By Corollary
\ref{c.i.grobner_extras}, $\Mult_{e_\gb}X_\ga^\gc$ is the number
of square-free monomials of maximal degree in $K[\cO_\gb]\setminus
\init G\abc$. By Lemma \ref{l.r.admiss_index}, this equals the
number of subsets $U$ of $\BBZ$ which are of maximal degree among
those bounded by $T_\ga, W_\gc$. However, a subset of $\BBZ$ is
bounded by $T_\ga, W_\gc$ if and only if it is bounded by
$\wt{T_\ga}, \wt{W_\gc}$ if and only if it is chain-bounded by
$\wt{T_\ga}, \wt{W_\gc}$, where the last equivalence is due to
Lemma \ref{l.t.two_orders}.
\end{proof}

\section{Path Families and Multiplicities}\label{s.path_families_mults}


For this section, we let $R$ and $S$ be fixed negative and
positive twisted chains contained in $\BBZ$ respectively. Let
\begin{align*}
\mathcal{M}_R&=\max\{U\subset \BBN\mid R\leqCC U\}\notag\\
\cM^S&=\max\{V\subset \BBP\mid V\leqCC S\}\notag\\
\cM_R^S&=\max\{W\subset \BBZ\mid R\leqCC W^- \text{ and }
W^+\leqCC S\}.
\end{align*}
where in each case by `max' we mean the subsets $U$, $V$, or $W$
respectively of maximal degree. For example, $\cM_R^S$ consists of
the collection of subsets $W$ of $\BBZ$ which are of maximal
degree among those which are chain-bounded by $R, S$. When
$R=\wt{T}_\ga$ and $S=\wt{W}_\gc$, $\cM_R^S$ consists precisely of
the subsets $U$ of Lemma \ref{l.chain_order}. In this section, in
order to give a better formulation of Lemma \ref{l.chain_order}
(see Proposition \ref{p.mult_chain}), we study the combinatorics
of $\cM_R^S$. Many of the definitions and ideas in this section
are illustrated in Examples \ref{e.t.multiplicities},
\ref{e.t.multiplicities2}, and \ref{e.t.multiplicities3}.

Note that
\begin{equation}\label{e.MRS}
\cM_R^S=\{U\dot{\cup} V\mid U\in \cM_R, V\in \cM^S\}.
\end{equation}
To study $\cM_R^S$, we begin by considering $\cM_R$, and thus
restricting attention to negative subsets of $\BBZ$. We say that a
subset $P\subset\BBN$ is \textbf{depth-one} if it contains no
two-element chains. If $P$ is depth-one, then we say that it is a
\textbf{negative-path} if the consecutive points are `as close as
possible' to each other, so that the points form a contiguous path
on $\BBN$ which moves only down or to the right.

For $r=(e,f)\in\BBN$, define
\begin{align*}
\floor{r} &= (e,f'),\text{ where }f'=\min\{y\in
\gb\mid (e,y)\in \BBN\}\\
\ceil{r} &= (e',f),\text{ where }e'=\max\{x\in \ol{\gb}\mid
(x,f)\in \BBN\}
\end{align*} We form the path $P_r$, which begins at $\floor{r}$,
moves horizontally to $r$, then vertically to $\ceil{r}$. Note
that since $R$ is a twisted chain,  if $r'\neq r$ then $P_{r'}\cap
P_{r}=\es$. Furthermore, $R\,\leqCC\, \dot{\bigcup}_{r\in R}P_r$.
Define $d_R=\sum_{r\in R} |P_r|$. The following lemma is a
straightforward consequence of the definitions.
\begin{lem}\label{l.paths_max}
Let $Q$ be a depth-one negative subset of $\BBZ$ such that
$P_r\leqCC Q$. Then $|Q|\leq |P_r|$, with equality if and only if
$Q$ is a negative-path from $\floor{r}$ to $\ceil{r}$.
\end{lem}

If $U\subset\BBN$, $R\leqCC U$, and $r\in R$, then define
$$U_{R,r}:=\{u\in U\mid r\leqCC u, \depth_U(u)=\depth_R(r)\}.$$
It follows from this definition that $U_{R,r}$ is depth-one.
Indeed, if $u$ and $u'$ are two elements of $U$ which form a
chain, then without loss of generality $u\ltC u'$. Thus
$\depth_U(u)<\depth_U(u')$, and in particular
$\depth_U(u)\neq\depth_U(u')$. Thus $u$ and $u'$ cannot both lie
in $U_{R,r}$.

\begin{lem}\label{l.paths_U} Let $U\subset \BBN$.

\noindent (i) If $R\leqCC U$, then $U=\dot{\bigcup}_{r\in
R}U_{R,r}$.

\noindent (ii) If $R\leqCC U$, then $|U|\leq d_R$, with equality
if and only if $U=\dot{\bigcup}_{r\in R}Q_r$, where $Q_r$ is a
negative-path from $\floor{r}$ to $\ceil{r}$.

\noindent (iii) Let $U=\dot{\bigcup}_{r\in R}Q_r\subset \BBN$,
where $Q_r$ is a negative-path from $\floor{r}$ to $\ceil{r}$.
Then $R\leqCC U$.

\end{lem}
\begin{proof}
(i) Let $u\in U$.  Then since $R\leqCC U$,
$\depth_U(u)\leq\depth_R(u)$.  Thus $\depth_U(u)=\depth_R(r)$ for
some $r\leqCC u$. This proves that $U=\bigcup_{r\in R}U_{R,r}$. To
prove that the union is disjoint, let $r,r'\in R$, $r\neq r'$, and
let $v\in U_{R,r}\cap U_{R,r'}$. By definition of $U_{R,r}$ and
$U_{R,r'}$, $\depth_R(r)=\depth_U(v)=\depth_R(r')$. Thus
$r\not\ltC r'$ and $r'\not\ltC r$. Since $R$ is a twisted chain,
$r\meet r'\not\in (\NN)^-$. But this implies that
$v\not\in(\NN)^-$, a contradiction.\\

\noindent (ii) For each $u\in U_{R,r}$, $r\leqCC u$; thus since
$U_{R,r}$ is depth-one, $P_r\leqCC U_{R,r}$. By Lemmas
\ref{l.paths_max} and \ref{l.paths_U}(i),
$|U|=|\dot{\bigcup}_{r\in R}U_{R,r}|=\sum_{r\in R}|U_{R,r}|\leq
\sum_{r\in R}|P_r|=d_R$, with equality if and only if for all
$r\in R$, $U_{R,r}$ is a negative-path from $\floor{r}$ to
$\ceil{r}$.  We denote $U_{R,r}$ by $Q_r$.\\

\noindent (iii) For each $r\in R$, $\{r\}\,\leqCC\, Q_r$.  Thus
$R\,\leqCC\, \dot{\bigcup}_{r\in R}r\,\leqCC\, \dot{\bigcup}_{r\in
R}Q_r=U$.
\end{proof}

Lemma \ref{l.paths_U}(ii) implies that any $U\in\cM_R$ is a
disjoint union $U=\dot{\bigcup}_{r\in R}Q_r$, where $Q_r$ is a
negative-path from $\floor{r}$ to $\ceil{r}$. Lemma
\ref{l.paths_U}(iii) implies that any disjoint union
$U=\dot{\bigcup}_{r\in R}Q_r$, where $Q_r$ is a negative-path from
$\floor{r}$ to $\ceil{r}$, is an element of $\cM_r$. Consequently
we have

\begin{cor}
$\cM_R$ consists of the set of all possible disjoint unions
$U=\dot{\bigcup}_{r\in R}Q_r$, where $Q_r$ is a negative-path from
$\floor{r}$ to $\ceil{r}$.
\end{cor}

Similar analysis can be done on positive subsets of $\BBZ$.  Here
the notion of a \textbf{positive-path} is identical to that of a
negative-path, except that it is contained in $\BBP$ instead of
$\BBN$. Likewise, the notions of $\floor{s}$ and $\ceil{s}$ for
$s\in\BBP$ are defined analogously as for $s\in\BBN$ (see Example
\ref{e.t.multiplicities}). We obtain

\begin{cor}
$\cM^S$ consists of the set of all possible disjoint unions
$V=\dot{\bigcup}_{s\in S}Q_s$, where $Q_s$ is a positive-path from
$\floor{s}$ to $\ceil{s}$.
\end{cor}
\noindent The preceding two corollaries and (\ref{e.MRS}) imply
\begin{cor}\label{c.path_families} $\cM_R^S$
consists of the set of all possible disjoint unions
$W=\dot{\bigcup}_{r\in R\cup S}Q_r$, where $Q_r$ is either a
negative-path or a positive-path from $\floor{r}$ to $\ceil{r}$,
depending on whether $r$ is negative or positive.
\end{cor}
\noindent The subsets $U$ of Lemma \ref{l.chain_order} are
precisely the elements of $\cM_R^S$, when $R=\wt{T}_\ga$ and
$S=\wt{W}_\gc$. Therefore combining Lemma \ref{l.chain_order} and
Corollary \ref{c.path_families}, we obtain
\begin{prop}\label{p.mult_chain}
$\Mult_{e_\gb}X_\ga^\gc$ is the number of disjoint unions
$\dot{\bigcup}_{r\in \wt{T_\ga}\cup \wt{W_\gb}}P_r$, where $P_r$
is either a negative-path or a positive-path from $\floor{r}$ to
$\ceil{r}$, depending on whether $r$ is negative or positive.
\end{prop}
\noindent We call a disjoint union as in Proposition
\ref{p.mult_chain} a family of nonintersecting paths in $\BBZ$.
\noindent Proposition \ref{p.mult_chain} also appears in
\cite{Ko-Ra}, \cite{Kra1}, \cite{Kra2}, and \cite{Kr-La}.

\eject

\begin{ex}\label{e.t.multiplicities}
Let $d=8$, $n=17$, $\ga=\{1,2,3,5,6,8,11,14\}$,\\
$\gb=\{2,7,8,9,12,13,16,17\}$, $\gc=\{8,9,11,13,14,15,16,17\}$.\\
(a) The negative and positive twisted chains
$\wt{T_\ga}=\{r_1,\ldots,r_6\}$ and
$\wt{W_\gc}=\{s_1,\ldots,s_3\}$ in $\gB$.\\
(b) The set of $\floor{r}$'s and $\ceil{r}$'s, for all
$r\in\wt{T_\ga}\cup \wt{W_{\gc}}$. Note that
$\floor{r_4}=\ceil{r_4}=r_4$ and
$\floor{r_5}=\ceil{r_5}=r_5$.\\
(c), (d) Two families of nonintersecting paths from $\floor{r}$ to
$\ceil{r}$, $r\in \wt{T_\ga}\cup \wt{W_{\gc}}$.
$\Mult_{e_\gb}X_\ga^\gc$ is the number of such families.  Note
that the path family in (c) consists of the paths $\{P_r\mid r\in
\wt{T_\ga}\cup \wt{W_{\gc}}\}$.
\end{ex}

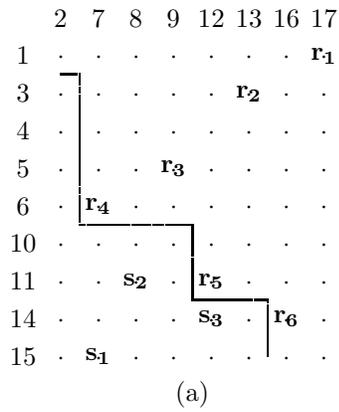
\begin{center}
\setlength{\unitlength}{.5mm}
$\begin{array}{l@{\hspace{1.5cm}}l}\\
%
\begin{picture}(90,100)(0,0)
\matrixput(10,10)(10,0){8}(0,10){9}{\plotsmcirchar}
\multiputlist(0,10)(0,10){15,14,11,10,6,5,4,3,1}
\multiputlist(10,100)(10,0){2,7,8,9,12,13,16,17}

\linethickness{.5pt}
\dottedline{.7}(10,85)(15,85)(15,45)(45,45)(45,25)(65,25)(65,10)

\put(80,90){\makebox(0,0){$\mathbf{r_1}$}}
\put(70,20){\makebox(0,0){$\mathbf{r_6}$}}
\put(60,80){\makebox(0,0){$\mathbf{r_2}$}}
\put(50,30){\makebox(0,0){$\mathbf{r_5}$}}
\put(40,60){\makebox(0,0){$\mathbf{r_3}$}}
\put(20,50){\makebox(0,0){$\mathbf{r_4}$}}

\put(20,10){\makebox(0,0){$\mathbf{s_1}$}}
\put(30,30){\makebox(0,0){$\mathbf{s_2}$}}
\put(50,20){\makebox(0,0){$\mathbf{s_3}$}}

\put(45,0){\makebox(0,0){\textrm{(a)}}}
\end{picture}
&

%
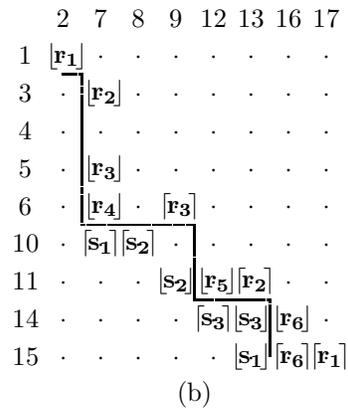
\begin{picture}(90,100)(0,0)
\matrixput(10,10)(10,0){8}(0,10){9}{\plotsmcirchar}
\multiputlist(0,10)(0,10){15,14,11,10,6,5,4,3,1}
\multiputlist(10,100)(10,0){2,7,8,9,12,13,16,17}

\linethickness{.5pt}
\dottedline{.7}(10,85)(15,85)(15,45)(45,45)(45,25)(65,25)(65,10)

\put(10,90){\makebox(2,0){$\mathbf{\sfloor{r_1}}$}}
\put(80,10){\makebox(2,0){$\mathbf{\sceil{r_1}}$}}

\put(70,20){\makebox(2,0){$\mathbf{\sfloor{r_6}}$}}
\put(70,10){\makebox(2,0){$\mathbf{\sceil{r_6}}$}}

\put(20,80){\makebox(2,0){$\mathbf{\sfloor{r_2}}$}}
\put(60,30){\makebox(2,0){$\mathbf{\sceil{r_2}}$}}

\put(50,30){\makebox(2,0){$\mathbf{\sfloor{r_5}}$}}

\put(20,60){\makebox(2,0){$\mathbf{\sfloor{r_3}}$}}
\put(40,50){\makebox(2,0){$\mathbf{\sceil{r_3}}$}}

\put(20,50){\makebox(2,0){$\mathbf{\sfloor{r_4}}$}}

\put(20,40){\makebox(0,0){$\mathbf{\sceil{s_1}}$}}
\put(60,10){\makebox(0,0){$\mathbf{\sfloor{s_1}}$}}

\put(30,40){\makebox(0,0){$\mathbf{\sceil{s_2}}$}}
\put(40,30){\makebox(0,0){$\mathbf{\sfloor{s_2}}$}}

\put(50,20){\makebox(0,0){$\mathbf{\sceil{s_3}}$}}
\put(60,20){\makebox(0,0){$\mathbf{\sfloor{s_3}}$}}

\put(45,0){\makebox(0,0){\textrm{(b)}}}
\end{picture}
\\[4em]

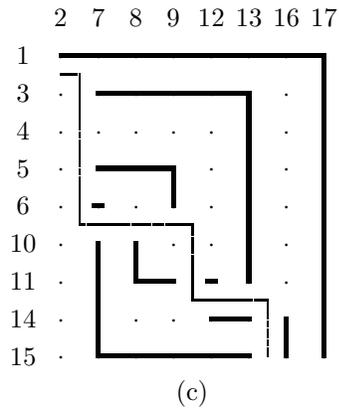
\begin{picture}(90,100)(0,0)
\matrixput(10,10)(10,0){8}(0,10){9}{\plotsmcirchar}
\multiputlist(0,10)(0,10){15,14,11,10,6,5,4,3,1}
\multiputlist(10,100)(10,0){2,7,8,9,12,13,16,17}

\linethickness{.5pt}
\dottedline{.7}(10,85)(15,85)(15,45)(45,45)(45,25)(65,25)(65,10)

%

\linethickness{1.5pt}
\dottedline{.7}(10,90)(80,90)(80,10)
\dottedline{.7}(20,80)(60,80)(60,30)
\dottedline{.7}(20,60)(40,60)(40,50) \dottedline{.7}(70,20)(70,10)

\dottedline{.7}(19,50)(21,50) \dottedline{.7}(49,30)(51,30)

\dottedline{.7}(20,40)(20,10)(60,10)
\dottedline{.7}(30,40)(30,30)(40,30) \dottedline{.7}(50,20)(60,20)

\put(45,0){\makebox(0,0){\textrm{(c)}}}
\end{picture}
&

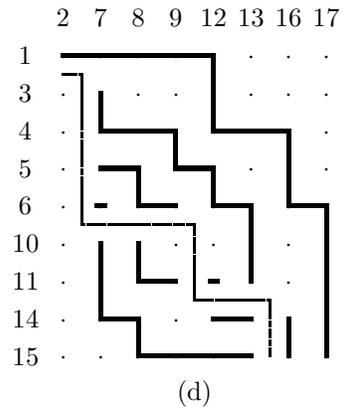
\begin{picture}(90,100)(0,0)
\matrixput(10,10)(10,0){8}(0,10){9}{\plotsmcirchar}
\multiputlist(0,10)(0,10){15,14,11,10,6,5,4,3,1}
\multiputlist(10,100)(10,0){2,7,8,9,12,13,16,17}

\linethickness{.5pt}
\dottedline{.7}(10,85)(15,85)(15,45)(45,45)(45,25)(65,25)(65,10)

%

\linethickness{1.5pt}
\dottedline{.7}(10,90)(50,90)(50,70)(70,70)(70,50)(80,50)(80,10)
\dottedline{.7}(20,80)(20,70)(40,70)(40,60)(50,60)(50,50)(60,50)(60,30)
\dottedline{.7}(20,60)(30,60)(30,50)(40,50)
\dottedline{.7}(70,20)(70,10)

\dottedline{.7}(19,50)(21,50)
\dottedline{.7}(49,30)(51,30)

\dottedline{.7}(20,40)(20,20)(30,20)(30,10)(60,10)
\dottedline{.7}(30,40)(30,30)(40,30)
\dottedline{.7}(50,20)(60,20)

\put(45,0){\makebox(0,0){\textrm{(d)}}}
\end{picture}

\end{array}
$
\end{center}

\eject

\begin{ex}\label{e.t.multiplicities2}
Let $d=4$, $n=9$, $\ga=\{1,2,3,5\}$, $\gb=\{1,5,6,8\}$,
$\gc=\{3,6,8,9\}$. We compute $\Mult_{e_\gb}X_\ga^\gc$. The
following two diagrams show the negative and positive twisted
chains $\wt{T_\ga}=\{r_1,r_2\}$ and $\wt{W_\gc}=\{s_1,s_2\}$ in
$\gB$; and the set of $\floor{r}$'s and $\ceil{r}$'s, for all
$r\in\wt{T_\ga}\cup \wt{W_{\gc}}$.

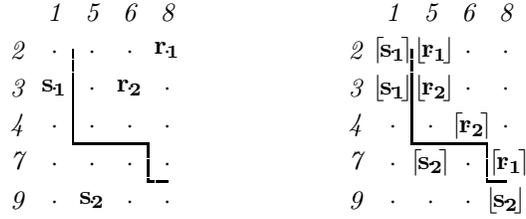
\setlength{\unitlength}{.5mm}
\begin{center}
$\begin{array}{l@{\hspace{2cm}}l}\\
\begin{picture}(50,60)(0,0)
\matrixput(10,10)(10,0){4}(0,10){5}{\plotsmcirchar}
\multiputlist(0,10)(0,10){9,7,4,3,2}
\multiputlist(10,60)(10,0){1,5,6,8}

\linethickness{.5pt}
\dottedline{.7}(15,50)(15,25)(35,25)(35,15)(40,15)

\put(40,50){\makebox(0,0){$\mathbf{r_1}$}}
\put(30,40){\makebox(0,0){$\mathbf{r_2}$}}
\put(10,40){\makebox(0,0){$\mathbf{s_1}$}}
\put(20,10){\makebox(0,0){$\mathbf{s_2}$}}
\end{picture}
&

\begin{picture}(50,60)(0,0)
\matrixput(10,10)(10,0){4}(0,10){5}{\plotsmcirchar}
\multiputlist(0,10)(0,10){9,7,4,3,2}
\multiputlist(10,60)(10,0){1,5,6,8}

\linethickness{.5pt}
\dottedline{.7}(15,50)(15,25)(35,25)(35,15)(40,15)

\put(20,50){\makebox(2,0){$\mathbf{\sfloor{r_1}}$}}
\put(40,20){\makebox(2,0){$\mathbf{\sceil{r_1}}$}}

\put(20,40){\makebox(2,0){$\mathbf{\sfloor{r_2}}$}}
\put(30,30){\makebox(2,0){$\mathbf{\sceil{r_2}}$}}

\put(10,50){\makebox(0,0){$\mathbf{\sceil{s_1}}$}}
\put(10,40){\makebox(0,0){$\mathbf{\sfloor{s_1}}$}}

\put(20,20){\makebox(0,0){$\mathbf{\sceil{s_2}}$}}
\put(40,10){\makebox(0,0){$\mathbf{\sfloor{s_2}}$}}
\end{picture}
\end{array}
$
\end{center}

\noindent There are six nonintersecting path families from
$\floor{r}$ to $\ceil{r}$, $r\in \wt{T_\ga}\cup \wt{W_{\gc}}$, as
shown below. Thus $\Mult_{e_\gb}X_\ga^\gc=6$.

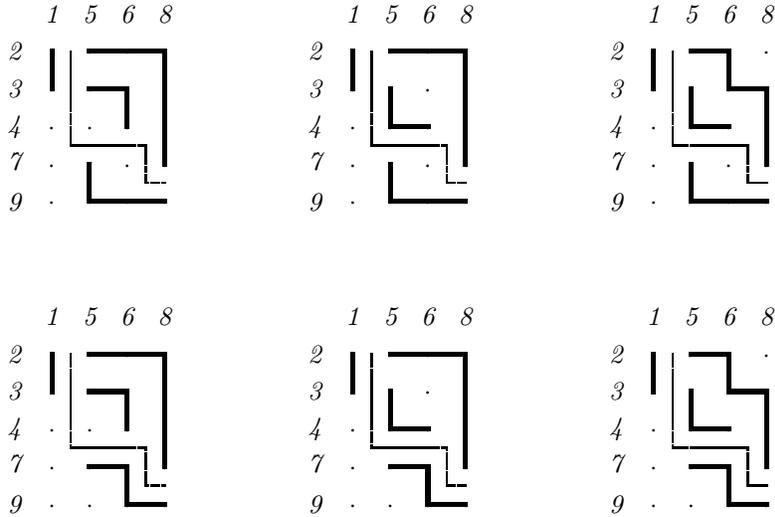
\begin{center}
$\begin{array}{l@{\hspace{1.5cm}}l@{\hspace{1.5cm}}l}\\
\begin{picture}(50,60)(0,0)
\matrixput(10,10)(10,0){4}(0,10){5}{\plotsmcirchar}
\multiputlist(0,10)(0,10){9,7,4,3,2}
\multiputlist(10,60)(10,0){1,5,6,8}

\linethickness{.5pt}
\dottedline{.7}(15,50)(15,25)(35,25)(35,15)(40,15)

\linethickness{1.5pt}
\dottedline{.7}(20,50)(40,50)(40,20)
\dottedline{.7}(20,40)(30,40)(30,30)
\dottedline{.7}(10,50)(10,40)
\dottedline{.7}(20,20)(20,10)(40,10)
\end{picture}
&

\begin{picture}(50,60)(0,0)
\matrixput(10,10)(10,0){4}(0,10){5}{\plotsmcirchar}
\multiputlist(0,10)(0,10){9,7,4,3,2}
\multiputlist(10,60)(10,0){1,5,6,8}

\linethickness{.5pt}
\dottedline{.7}(15,50)(15,25)(35,25)(35,15)(40,15)

\linethickness{1.5pt}
\dottedline{.7}(20,50)(40,50)(40,20)
\dottedline{.7}(20,40)(20,30)(30,30)
\dottedline{.7}(10,50)(10,40)
\dottedline{.7}(20,20)(20,10)(40,10)
\end{picture}
&

\begin{picture}(50,60)(0,0)
\matrixput(10,10)(10,0){4}(0,10){5}{\plotsmcirchar}
\multiputlist(0,10)(0,10){9,7,4,3,2}
\multiputlist(10,60)(10,0){1,5,6,8}

\linethickness{.5pt}
\dottedline{.7}(15,50)(15,25)(35,25)(35,15)(40,15)

\linethickness{1.5pt}
\dottedline{.7}(20,50)(30,50)(30,40)(40,40)(40,20)
\dottedline{.7}(20,40)(20,30)(30,30)
\dottedline{.7}(10,50)(10,40)
\dottedline{.7}(20,20)(20,10)(40,10)
\end{picture}
\\[2.5em]

\begin{picture}(50,60)(0,0)
\matrixput(10,10)(10,0){4}(0,10){5}{\plotsmcirchar}
\multiputlist(0,10)(0,10){9,7,4,3,2}
\multiputlist(10,60)(10,0){1,5,6,8}

\linethickness{.5pt}
\dottedline{.7}(15,50)(15,25)(35,25)(35,15)(40,15)

\linethickness{1.5pt}
\dottedline{.7}(20,50)(40,50)(40,20)
\dottedline{.7}(20,40)(30,40)(30,30)
\dottedline{.7}(10,50)(10,40)
\dottedline{.7}(20,20)(30,20)(30,10)(40,10)
\end{picture}
&

\begin{picture}(50,60)(0,0)
\matrixput(10,10)(10,0){4}(0,10){5}{\plotsmcirchar}
\multiputlist(0,10)(0,10){9,7,4,3,2}
\multiputlist(10,60)(10,0){1,5,6,8}

\linethickness{.5pt}
\dottedline{.7}(15,50)(15,25)(35,25)(35,15)(40,15)

\linethickness{1.5pt}
\dottedline{.7}(20,50)(40,50)(40,20)
\dottedline{.7}(20,40)(20,30)(30,30)
\dottedline{.7}(10,50)(10,40)
\dottedline{.7}(20,20)(30,20)(30,10)(40,10)
\end{picture}
&

\begin{picture}(50,60)(0,0)
\matrixput(10,10)(10,0){4}(0,10){5}{\plotsmcirchar}
\multiputlist(0,10)(0,10){9,7,4,3,2}
\multiputlist(10,60)(10,0){1,5,6,8}

\linethickness{.5pt}
\dottedline{.7}(15,50)(15,25)(35,25)(35,15)(40,15)

\linethickness{1.5pt}
\dottedline{.7}(20,50)(30,50)(30,40)(40,40)(40,20)
\dottedline{.7}(20,40)(20,30)(30,30)
\dottedline{.7}(10,50)(10,40)
\dottedline{.7}(20,20)(30,20)(30,10)(40,10)
\end{picture}

\end{array}
$
\end{center}

\end{ex}

\eject

\begin{ex}\label{e.t.multiplicities3}
Let $d=4$, $n=9$, $\ga=\{1,2,3,5\}$, $\gb=\{1,5,6,8\}$,
$\gc=\{3,6,8,9\}$. We compute multiplicities at $e_\gb$ of the
Schubert variety $X^{\gc}$ and the opposite Schubert variety
$X_\ga$. Note that $\ga$, $\gb$, and $\gc$ are the same as in the
previous example.

Observe that $X^{\gc}=X^{\gc}_{\text{id}}$ and
$X_{\ga}=X^{\omega_0}_{\ga}$, where $\text{id}=\{1,2,3,4\}$ and
$\omega_0=\{6,7,8,9\}$. We have that
$\wt{T_\text{id}}=\{\gs_{2,8},\gs_{3,6},\gs_{4,5}\}$ and
$\wt{W_{\omega_0}}=\{\gs_{9,1},\gs_{7,5}\}$, where $\gs_{i,j}$ is
the transposition exchanging $i$ and $j$. (Both $\wt{T_\text{id}}$
and $\wt{W_{\omega_0}}$ are in fact chains.)

There are two nonintersecting path families from $\floor{r}$ to
$\ceil{r}$, $r\in \wt{T_{\text{id}}}\cup \wt{W_{\gc}}$, as shown
below. Thus $\Mult_{e_\gb}X^\gc=\Mult_{e_\gb}X_{\text{id}}^\gc=2$.

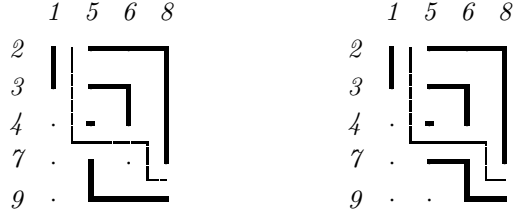
\setlength{\unitlength}{.5mm}
\begin{center}
$\begin{array}{l@{\hspace{2cm}}l}\\
\begin{picture}(50,60)(0,0)
\matrixput(10,10)(10,0){4}(0,10){5}{\plotsmcirchar}
\multiputlist(0,10)(0,10){9,7,4,3,2}
\multiputlist(10,60)(10,0){1,5,6,8}

\linethickness{.5pt}
\dottedline{.7}(15,50)(15,25)(35,25)(35,15)(40,15)

\linethickness{1.5pt}
\dottedline{.7}(20,50)(40,50)(40,20)
\dottedline{.7}(20,40)(30,40)(30,30)
\dottedline{.7}(19.5,30)(20.5,30)
\dottedline{.7}(10,50)(10,40) \dottedline{.7}(20,20)(20,10)(40,10)
\end{picture}
&

\begin{picture}(50,60)(0,0)
\matrixput(10,10)(10,0){4}(0,10){5}{\plotsmcirchar}
\multiputlist(0,10)(0,10){9,7,4,3,2}
\multiputlist(10,60)(10,0){1,5,6,8}

\linethickness{.5pt}
\dottedline{.7}(15,50)(15,25)(35,25)(35,15)(40,15)

\linethickness{1.5pt}
\dottedline{.7}(20,50)(40,50)(40,20)
\dottedline{.7}(20,40)(30,40)(30,30)
\dottedline{.7}(19.5,30)(20.5,30)
\dottedline{.7}(10,50)(10,40)
\dottedline{.7}(20,20)(30,20)(30,10)(40,10)
\end{picture}
\end{array}
$
\end{center}

\noindent There are three nonintersecting path families from
$\floor{r}$ to $\ceil{r}$, $r\in \wt{T_\ga}\cup
\wt{W_{\omega_0}}$, as shown below. Thus
$\Mult_{e_\gb}X_\ga=\Mult_{e_{\gb}}X_{\ga}^{\omega_0}=3$.

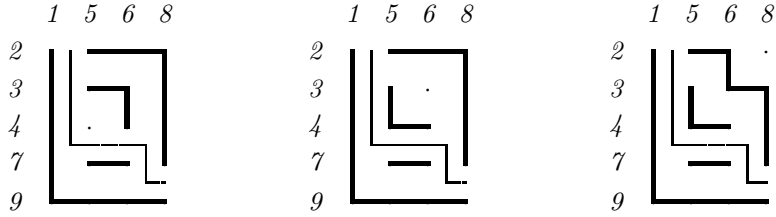
\begin{center}
$\begin{array}{l@{\hspace{1.5cm}}l@{\hspace{1.5cm}}l}\\
\begin{picture}(50,60)(0,0)
\matrixput(10,10)(10,0){4}(0,10){5}{\plotsmcirchar}
\multiputlist(0,10)(0,10){9,7,4,3,2}
\multiputlist(10,60)(10,0){1,5,6,8}

\linethickness{.5pt}
\dottedline{.7}(15,50)(15,25)(35,25)(35,15)(40,15)

\linethickness{1.5pt}
\dottedline{.7}(20,50)(40,50)(40,20)
\dottedline{.7}(20,40)(30,40)(30,30)
\dottedline{.7}(10,50)(10,10)(40,10)
\dottedline{.7}(20,20)(30,20)
\end{picture}
&

\begin{picture}(50,60)(0,0)
\matrixput(10,10)(10,0){4}(0,10){5}{\plotsmcirchar}
\multiputlist(0,10)(0,10){9,7,4,3,2}
\multiputlist(10,60)(10,0){1,5,6,8}

\linethickness{.5pt}
\dottedline{.7}(15,50)(15,25)(35,25)(35,15)(40,15)

\linethickness{1.5pt}
\dottedline{.7}(20,50)(40,50)(40,20)
\dottedline{.7}(20,40)(20,30)(30,30)
\dottedline{.7}(10,50)(10,10)(40,10)
\dottedline{.7}(20,20)(30,20)
\end{picture}
&

\begin{picture}(50,60)(0,0)
\matrixput(10,10)(10,0){4}(0,10){5}{\plotsmcirchar}
\multiputlist(0,10)(0,10){9,7,4,3,2}
\multiputlist(10,60)(10,0){1,5,6,8}

\linethickness{.5pt}
\dottedline{.7}(15,50)(15,25)(35,25)(35,15)(40,15)

\linethickness{1.5pt}
\dottedline{.7}(20,50)(30,50)(30,40)(40,40)(40,20)
\dottedline{.7}(20,40)(20,30)(30,30)
\dottedline{.7}(10,50)(10,10)(40,10)
\dottedline{.7}(20,20)(30,20)
\end{picture}
\end{array}
$
\end{center}

From this example and the previous one, we see that
$\Mult_{e_\gb}X_\ga^\gc$ equals the product
$\Mult_{e_\gb}X^\gc\cdot \Mult_{e_\gb}X_\ga$. It is not difficult
to show that this equality holds in general. This fact is also
proven in \cite{Kr-La3}, using different methods.
\end{ex}
\eject

\section{Proofs}\label{s.proofs}

In this section, we give proofs of Lemmas \ref{l.g.dominance} and
\ref{l.t.two_orders}.


\subsection*{Proof of Lemma \ref{l.g.dominance}}

Let $U$ be a nonvanishing multiset on $\NN$, and let $T$ and $W$
be negative and positive subsets of $\NN$ respectively with the
property that $T_{(1)},T_{(2)},W_{(1)}$, and $W_{(2)}$ are both
subsets of $\bN$, i.e., multisets such that each value has
cardinality one. Lemma \ref{l.g.dominance} is part (v) of the
following.
\begin{lem}
(i) Suppose that $U=\{(a_1,b_1),\ldots,(a_t,b_t)\}$ is a negative
multiset on $\NN$ whose entries are listed in lexicographic order.
For $k=1,\ldots,t$, let $U^{(k)}:=\{(a_1,b_1),\ldots,(a_k,b_k)\}$,
and let $(P^{(k)},Q^{(k)})=\rsk(U^{(k)})$ (note that
$(P^{(t)},Q^{(t)})=\rsk(U)$). Define
$\{p^{(k)}_1,\ldots,p^{(k)}_{c_k}\}$ to be the first row of
$P^{(k)}$ and $\{q^{(k)}_1,\ldots,q^{(k)}_{c_k}\}$ the first row
of $Q^{(k)}$, both listed in increasing order. Let
$m(k):=\max\{m\in\{1,\ldots,c_k\}\mid
p^{(k)}_{m}<q^{(k)}_1\}=|(P^{(k)}_1)^{<q^{(k)}_1}|$. Then for
$1\leq j\leq m(k)$, there exists a chain $C_{k,j}$ in $U^{(k)}$
which has
$j$ elements, the last of which has first component $p^{(k)}_j$.\\
\noindent (ii) If $U$ is bounded by $T,\es$, then
$(P^{(k)},Q^{(k)})$ is bounded by $T,\es$, $k=1,\ldots,t$.\\
\noindent (iii) If $U$ is
bounded by $T,\es$, then $\rsk(U)$ is bounded by $T,\es$.\\
\noindent (iv) If $U$ is bounded by $\es,W$, then $\rsk(U)$ is
bounded by $\es,W$.\\
\noindent (v) If $U$ is bounded by $T,W$, then $\rsk(U)$ is
bounded by $T,W$.
\end{lem}
\begin{proof}
We prove (i) and (ii) together by induction on $k$, with $k=1$ the
starting point for the induction. This case is technically covered
by Case 2 below, although it can be checked quite easily:
$U^{(1)}=\{(a_1,b_1)\}$, $P^{(1)}$ contains the sole entry $a_1$,
and $Q^{(1)}$ contains the sole entry $b_1$. For (i), $m(1)=1$,
$C_{1,1}=(a_1,b_1)$.  For (ii), $U$ bounded by $T,\es$ implies
$U^{(1)}$ bounded by $T,\es$, which is clearly equivalent to
$(P^{(1)},Q^{(1)})$ bounded by $T,\es$.

Let $k\in 2,\ldots,t-1$.  Let $(P,Q)=(P^{(k)},Q^{(k)})$,
$a=a_{k+1}$, $b=b_{k+1}$, $(P',Q')=(P^{(k+1)},Q^{(k+1)})$,
$U=U^{(k)}$, $U'=U^{(k+1)}$,
$\{p_1,\ldots,p_c\}=\{p^{(k)}_1,\ldots,p^{(k)}_{c_k}\}$,
$\{q_1,\ldots,q_c\}=\{q^{(k)}_1,\ldots,q^{(k)}_{c_k}\}$. Note that
$\{p_1,\ldots,p_c\}\subset \{a_1,\ldots,a_k\}$,
$\{q_1,\ldots,q_c\}\subset \{b_1,\ldots,b_k\}$. Thus since $b$ is
less than or equal to all elements of $\{b_1,\ldots,b_k\}$,
$a<b\leq q_1$. We assume inductively that $$T_{(1)}-T_{(2)}\leq
P_1-Q_1,$$ and prove that
$$T_{(1)}-T_{(2)}\leq P'_1-Q'_1.$$ Equivalently, we prove that
for all positive integers $z$, $$|(T_{(1)}-T_{(2)})^{\leq z}|
\geq|(P'_1-Q'_1)^{\leq z}|,$$ where we use the definition
$A-B:=A\disjunion (\bN\setminus B)$, where $A$ and $B$ are both
subsets of $\bN$ (see Section \ref{s.multisets}).

We consider two cases, corresponding to the two
ways in which $(P'_1,Q'_1)$ can be obtained from $(P_1,Q_1)$.\\

\noindent {\it Case 1}.  $P'_1$ is obtained by $a$ bumping $p_l$
in $P_1$, for some $1\leq l\leq c$, i.e.,
\[
\begin{array}{l}
P'_1=P_1\setminus \{p_l\}\disjunion \{a\}\\
Q'_1=Q_1
\end{array}
\]

\noindent (i) The fact that $a$ bumps $p_l$ implies both $a\leq
p_l$ and $p_l<b$.  Hence $a\leq p_l<b\leq q_1$, which implies
$(P'_1)^{<q_1}=(P_1)^{<q_1}\setminus \{p_l\}\cup \{a\}$.  Thus
$m(k+1)=m(k)$. For $j\in\{1,\ldots,m(k)\}\setminus\{l\}$, set
$C_{k+1,j}=C_{k,j}$. If $l=1$ then set $C_{k+1,l}=\{(a,b)\}$.
Otherwise, consider the chain
$C_{k,l-1}=\{(g_1,h_1),\ldots,(p_{l-1},h_{l-1})\}$.  Since $a$
bumps $p_l$, $a>p_{l-1}$ .  Thus $b<h_{l-1}$, since $(a,b)$ comes
after $(p_{l-1},h_{l-1})$ in the ordered list of elements of $U'$.
Therefore $C:=C_{k,l-1}\cup \{(a,b)\}$ is a chain
in $U'$.  We let $C_{k+1,l}$ be this chain.\\

\noindent (ii) For $z<a$ or $z\geq p_l$,
\begin{equation}\label{e.chd1}
|(T_{(1)}-T_{(2)})^{\leq z}|\geq|(P_1-Q_1)^{\leq
z}|=|(P'_1-Q'_1)^{\leq z}|.
\end{equation}
If $a=p_l$ then we are done.  Thus we assume that $a<p_l$. We
claim that for $a\leq z<p_l$, $ |(C_{(1)}-C_{(2)})^{\leq
z}|=|(P'_1-Q'_1)^{\leq z}|$.  Assuming the claim (and using the
fact that $T\leq C\leq\es$, since $U$ is bounded by $T,\es$) we
have that for $a\leq z<p_l$,
\begin{equation}\label{e.chd2}
|(T_{(1)}-T_{(2)})^{\leq z}|\geq|(C_{(1)}-C_{(2)})^{\leq
z}|=|(P'_1-Q'_1)^{\leq z}|.
\end{equation}
Now (\ref{e.chd1}) and (\ref{e.chd2}) prove the inductive step of
(ii).

We now prove the claim. From the proof of (i), we have that
$C=C_{k+1,l}$$=\{(g_1,h_1),\ldots,(g_{l-1},h_{l-1}),(a,b)\}$,
where $g_1<\cdots<g_{l-1}<a< p_l<b<h_{l-1}<\cdots<h_1$. Thus for
$a\leq z<p_l$,
$$|(C_{(1)}-C_{(2)})^{\leq z}|=|(C_{(1)}\disjunion (\bN\setminus
C_{(2)}))^{\leq z}|=|(C_{(1)}\disjunion \bN)^{\leq z}|=l+z.$$
Also, $p_1<\cdots<p_{l-1}<a< p_l<b\leq q_1<\cdots<q_c$. Thus for
$a\leq z<p_l$,
$$|(P'_1-Q'_1)^{\leq z}|=|(P'_1\disjunion (\bN\setminus Q'_1))^{\leq z}|=
|(P'_1\disjunion \bN)^{\leq z}|=l+z.$$

\vspace{1em}

\noindent {\it Case 2}.  $P'_1$ is obtained by adding $a$ to $P_1$
in position $l$ from the left and $Q'_1$ is obtained by adding $b$
to the left end of $Q_1$ (and shifting all other entries of $Q_1$
to the right by one box), i.e.,
\[
\begin{array}{l}
P'_1=P_1\disjunion \{a\}\\
Q'_1=Q_1\disjunion \{b\}
\end{array}
\]
\noindent (i) We have that
$P'_1=\{p_1,\ldots,p_{l-1},a,p_{l},\ldots,p_c\}$,
$Q'_1=\{b,q_1,\ldots,q_c\}$, where the elements of both sets are
listed in strictly increasing order (note that $b<q_1$ follows
from the fact that $Q'$ is row-strict, which is proven in Lemma
\ref{l.ins_preserves_semist_bitab}). Now $p_{l-1}<a<b<q_1$ implies
that $m(k)\geq l-1$.  Note that $a<b\leq p_l$ (since $b>p_l$ would
require that $a$ bump $p_l$ in the bounded insertion process).
Thus $m(k+1)=l$. For $j\in\{1,\ldots,l-1\}$, set
$C_{k+1,j}=C_{k,j}$. Consider the chain
$C_{k,l-1}=\{(g_1,h_1),\ldots,(g_{l-2},h_{l-2}),(p_{l-1},h_{l-1})\}$.
Now $a>p_{l-1}$, and this implies that $b<h_{l-1}$, since $(a,b)$
comes after $(p_{l-1},h_{l-1})$ in the ordered list of elements of
$U'$. Therefore $C:=C_{k,l-1}\cup \{(a,b)\}$ is a chain in $U'$.
We let $C_{k+1,l}$ be this chain.

\vspace{1em}

\noindent (ii) For $z<a$,
\begin{equation}\label{e.chd21}
|(T_{(1)}-T_{(2)})^{\leq z}|\geq|(P_1-Q_1)^{\leq
z}|=|(P'_1-Q'_1)^{\leq z}|.
\end{equation}
In fact, (\ref{e.chd21}) holds for $z\geq b$ as well, since for
such $z$,
\begin{align*}
|(P'_1-Q'_1)^{\leq z}|&=|(P'_1\disjunion (\bN\setminus
Q'_1))^{\leq z}|\\
&=|(P'_1)^{\leq z}|+|(\bN\setminus Q'_1)^{\leq
z}|\\
&=(|(P_1)^{\leq z}|+1)+(|(\bN\setminus Q_1)^{\leq
z}|-1)\\
&=|(P_1)^{\leq z}|+|(\bN\setminus Q_1)^{\leq z}|\\
&=|(P_1-Q_1)^{\leq z}|.
\end{align*}
We claim that for $a\leq z<b$, $ |(C_{(1)}-C_{(2)})^{\leq
z}|=|(P'_1-Q'_1)^{\leq z}|$.  Assuming the claim (and using the
fact that $T\leq C\leq\es$, since $U$ is bounded by $T,\es$) we
have that for $a\leq z<b$,
\begin{equation}\label{e.chd22}
|(T_{(1)}-T_{(2)})^{\leq z}|\geq|(C_{(1)}-C_{(2)})^{\leq
z}|=|(P'_1-Q'_1)^{\leq z}|.
\end{equation}
Now (\ref{e.chd21}) and (\ref{e.chd22}) prove the inductive step
of (ii).

We now prove the claim. From the proof of (i), we have that
$C=C_{k+1,l}$$=\{(g_1,h_1),\ldots,(g_{l-1},h_{l-1}),(a,b)\}$,
where $g_1<\cdots<g_{l-1}<a<b<h_{l-1}<\cdots<h_1$. Thus for $a\leq
z<b$,
$$|(C_{(1)}-C_{(2)})^{\leq z}|=|(C_{(1)}\disjunion (\bN\setminus
C_{(2)}))^{\leq z}|=|(C_{(1)}\disjunion \bN)^{\leq z}|=l+z.$$
Also, $p_1<\cdots<p_{l-1}<a<b\leq  p_l$, $b< q_1<\cdots<q_c$. Thus
for $a\leq z<b$,
$$|(P'_1-Q'_1)^{\leq z}|=|(P'_1\disjunion (\bN\setminus Q'_1))^{\leq z}|=
|(P'_1\disjunion \bN)^{\leq z}|=l+z.$$

\noindent (iii) Set $k=t$ in (ii).

\noindent (iv) Use arguments similar to (i), (ii), and (iii), but
for $U$ positive. Alternatively, one could apply the involution
$\iota$ to (iii).

\noindent (v) Use (iii), (iv), and the fact that $U$ is bounded by
$T,W$ if and only if $U^-$ is bounded by $T,\es$ and $U^+$ is
bounded by $\es,W$; and similarly for $\rsk(U)$.

\end{proof}


\subsection*{Proof of Lemma \ref{l.t.two_orders}}

Parts (iii) and (iv) of the Lemma below imply Lemma
\ref{l.t.two_orders}. In this proof, for $R$ a subset of $\NN$, we
define $R_{(1)}-R_{(2)}$ to be the (infinite) multiset
$R_{(1)}\disjunion\bN\setminus R_{(2)}$ (see Section
\ref{s.multisets}).
\begin{lem}
(i) Let $R, S$ be negative twisted chains. Then $R\leqCC S$ if and
only if $\depth_R((z,z+1))\geq \depth_S((z,z+1))$ for all
$z\in\bN$.\\
(ii) Let $R$ be a negative twisted chain, and let $z\in\bN$. Then
$\depth_R(z,z+1)=|(R_{(1)}-R_{(2)})^{\leq
z}|-z$.\\
(iii) Let $R,S$ be negative twisted chains. Then
$R\leqCC S\iff R\leq S$.\\
(iv) Let $R,S$ be positive twisted chains.  Then
$R\leqCC S\iff R\leq S$.\\
\end{lem}
\begin{proof}
(i) The ``only if'' direction is obvious. For $(e,f)\in(\bN^2)^-$,
define $D((e,f))=\depth_R((e,f))-\depth_S((e,f))$.  Suppose that
$D((e,f))<0$ for some $(e,f)\in(\bN^2)^-$. We must show that
$D((z,z+1))<0$ for some $z\in \bN$.  Suppose that $D((e,e+1))\geq
0$. Then $\depth_R((e,e+1))>\depth_R((e,f))$. Thus, there exists
$(e',f')\in R$ such that $e'\leq e$ and $e+1\leq f'< f$. Let
$(g,h)$ be the one such with maximal $f'$.

We claim that $\depth_R((h,h+1))=\depth_R((e,f))$.  If not, then
there exists $(p,q)\in R$ such that either (a) $e<p< h$, $h+1\leq
q$, (b) $p=h$, $h+1\leq q$, or (c) $p\leq e$, $h+1\leq q<f$. In
case (a), $(p,q)\not\ltC (g,h)$, $(g,h)\not\ltC (p,q)$, and
$(g,h)\wedge(p,q)=(p,h)\in(\NN)^-$ (since $p<h$), contradicting
the fact that $R$ is a negative twisted chain. In case (b), $p=h$,
and thus $(g,h),(h,q)\in R$ contradicts the fact that $R$ is
completely disjointed. In case (c), the maximality of $h$ is
violated.

Since $(e,f)\leqCC(h,h+1)$,
$\depth_S((h,h+1))\geq\depth_S((e,f))$. Thus $D((h,h+1))\leq
D((e,f))<0$.\\

\noindent (ii) Let $R=\{(e_1,f_1),\ldots,(e_m,f_m)\}$, with
$e_1<\cdots<e_m$. Note that $\{(e_i,f_i)\in R\mid e_i\leq
z<z+1\leq f_i\}$ consists of all the $(e_i,f_i)\in R$ such that
$(e_i,f_i)\leqCC (z,z+1)$. Thus, since $R$ is a twisted chain,
$\{(e_i,f_i)\in R\mid e_i\leq z<z+1\leq f_i\}$ must form a chain.
Hence $\depth_R(z,z+1)=|\{(e_i,f_i)\in R\mid e_i\leq z<z+1\leq
f_i\}|$.

Recall that
$R_{(1)}-R_{(2)}=\{e_1,\ldots,e_m\}\disjunion(\bN\setminus
\{f_1,\ldots,f_m\})$. The result now follows from the fact that
$$|\{(e_i,f_i)\in R\mid e_i\leq z<z+1\leq
f_i\}|=|(R_{(1)}-R_{(2)})^{\leq z}\setminus \bN^{\leq z}|.$$ To
see why this equality holds, observe that $R_{(1)}-R_{(2)}$ can be
obtained by starting with $\bN$ and then successively replacing
$f_i$ by $e_i$, $i=1,\ldots,m$. Such a replacement adds 1 to the
number of elements less than or equal to $z$ if and only if
$e_i\leq z$ and $f_i\geq z+1$ (and never subtracts 1 from the
number of elements less than or equal to $z$, since $e_i<f_i$).\\

\noindent (iii)
\begin{align*}
R\leqCC S&\iff \depth_R((z,z+1))\geq\depth_S((z,z+1)),\ z\in\bN\\
&\iff |(R_{(1)}-R_{(2)})^{\leq z}|\geq |(S_{(1)}-S_{(2)})^{\leq z}|,\ z\in\bN \\
&\iff R_{(1)}-R_{(2)}\leq S_{(1)}-S_{(2)}\\
&\iff R\leqT S,
\end{align*}
where the first and second equivalences follow from (i) and
(ii) respectively.\\

\noindent (iv) is obtained by applying $\iota$ to both sides of
(iii).

\end{proof}


\providecommand{\bysame}{\leavevmode\hbox
to3em{\hrulefill}\thinspace}
\providecommand{\MR}{\relax\ifhmode\unskip\space\fi MR }
\providecommand{\MRhref}[2]{%
  \href{http://www.ams.org/mathscinet-getitem?mr=#1}{#2}
} \providecommand{\href}[2]{#2}

\vspace{1em}

\noindent \textsc{Department of Mathematics, Univ. of Georgia,
Athens, GA 30602}

\noindent \textsl{Email address}: \texttt{vkreiman@math.uga.edu}

\noindent May 18, 2007
\end{document}